\documentclass[11pt]{article}

\usepackage[margin=1in]{geometry}
\usepackage{amsmath,amssymb,amsthm,mathtools}
\usepackage{microtype}
\usepackage{booktabs}
\usepackage{subcaption}
\usepackage{graphicx}
\usepackage[inline]{enumitem}
\usepackage{multirow,longtable,tabularx,makecell}
\usepackage[flushleft]{threeparttable}
\usepackage{bm,bbm,bigints,braket}
\usepackage{algorithm}
\usepackage{algpseudocode}
\usepackage{tikz}
\usepackage{pgfplots}
\usepackage{thmtools}
\usepackage{thm-restate}
\usetikzlibrary{arrows.meta,backgrounds,calc,decorations.pathreplacing,patterns,plotmarks,positioning}
\usepgflibrary{arrows}
\makeatletter
\def\theHALG@line{\arabic{algorithm}.\arabic{ALG@line}}
\makeatother
\usepackage{natbib}
\usepackage[hidelinks]{hyperref}
\usepackage[nameinlink,capitalize,noabbrev]{cleveref}
\declaretheorem[name=Theorem,style=plain]{theorem}
\declaretheorem[name=Proposition,style=plain]{proposition}
\declaretheorem[name=Lemma,style=plain]{lemma}

\declaretheorem[name=Corollary,style=plain]{corollary}
\declaretheorem[name=Assumption,style=plain]{assumption}

\declaretheorem[name=Example,style=definition]{example}

\def\E{{\mathbb E}}

\def\Pr{{\mathbb{P}}}

\def\Re{\mathbb{R}}

\def\I{\mathbb{I}}

\def\hat{\widehat}

\def \P{\mathcal{P}}

\def \Ze{{\mathbb{Z}}}

\def\O{{\mathcal O}}

\def\N{{\mathcal N}}

\def\P{{\mathcal P}}

\def\X{{\mathcal X}}

\def\Re{{\mathbb R}}

\def\O{{\mathcal O}}

\def\CVaR{{\mathrm{CVaR}}}
\def\VaR{{\mathrm{VaR}}}

\def\as{A\#}

\DeclarePairedDelimiter\ceil{\lceil}{\rceil}

\newcommand{\trxi}{\tilde{\bm{\xi}}}
\newcommand{\rxi}{\bm{\xi}}

\pgfplotsset{compat=newest}
\pgfdeclareplotmark{mystar}{
\node[star,star point ratio=2.25,minimum size=6pt,
inner sep=0pt,draw=black,solid,fill=red] {};
}
\algdef{SE}[DOWHILE]{Do}{doWhile}{\algorithmicdo}[1]{\algorithmicwhile\ #1}
\algnewcommand{\Or}{\textbf{or}}
\algnewcommand{\And}{\textbf{and}}
\DeclareMathOperator*{\argmin}{arg\,min}
\def\x{\vect{x}}
\newcommand{\vect}[1]{\boldsymbol{\bm{#1}}}
\newcommand*{\QEDA}{\qedhere}
\newcommand*{\QEDB}{\hfill\ensuremath{\diamond}}
\newenvironment{APPENDICES}{\appendix}{}

\definecolor{bluegreen}{RGB}{0,158,115}

\bibpunct[,]{(}{)}{,}{a}{}{,}

\allowdisplaybreaks

\title{On Sequential Convex Approximations of Wasserstein Distributionally Robust Chance Constrained Programs}
\author{
  \begin{tabular}{@{}ccc@{}}
  Yannan Wang\textsuperscript{1} & Nan Jiang\textsuperscript{1} & Haoming Shen\textsuperscript{2} \\
  \small\href{mailto:ywangwa@connect.ust.hk}{\nolinkurl{ywangwa@connect.ust.hk}} & \small\href{mailto:nanjiang@ust.hk}{\nolinkurl{nanjiang@ust.hk}} & \small\href{mailto:haomings@uark.edu}{\nolinkurl{haomings@uark.edu}}
  \end{tabular} \\[0.6em]
  \small \textsuperscript{1}Department of Industrial Engineering and Decision Analytics, \\
  \small The Hong Kong University of Science and Technology \\[0.2em]
  \small \textsuperscript{2}Department of Industrial Engineering, University of Arkansas
}
\date{}
\hypersetup{
  pdftitle={On Sequential Convex Approximations of Wasserstein Distributionally Robust Chance Constrained Programs},
  pdfauthor={Yannan Wang, Nan Jiang, Haoming Shen},
  pdfsubject={Sequential Convex Approximations of Wasserstein DRCCP},
  pdfkeywords={Distributionally Robust, Chance Constraints, Wasserstein, Sequential Convex Approximation, CVaR}
}

\begin{document}
\markboth{Wang, Jiang, and Shen}{Sequential Convex Approximations of Wasserstein DRCCP}
\maketitle
\begin{abstract}
We study distributionally robust chance constrained programs (DRCCPs) under the Wasserstein ambiguity set. Exact reformulations of such problems are generally nonconvex and NP-hard, while commonly used CVaR-based convex approximations can be overly conservative. To address these challenges, we derive an equivalent difference-of-convex reformulation of Wasserstein DRCCPs and exploit this structure to develop a sequential convex approximation algorithm. The proposed algorithm generates feasible solutions with monotone objective improvement, converges to critical points, and admits an explicit iteration-complexity guarantee. We further show that, under mild conditions, the algorithm provably improves upon the CVaR approximation after very few iterations. 
Numerical experiments show that the proposed method substantially improves both solution quality and running time compared to benchmark approaches.
\par\medskip
\noindent\textbf{Keywords:} Distributionally Robust, Chance Constraints, Wasserstein, Sequential Convex Approximation, CVaR
\end{abstract}

\section{Introduction}\label{sec_intro}
Distributionally robust chance constrained programming (DRCCP) has emerged as a powerful framework for decision-making under uncertainty, particularly when the underlying probability distribution is only partially known. This framework bridges regular chance constrained programming and distributionally robust optimization by simultaneously enforcing probabilistic feasibility and robustness against distributional ambiguity. Accordingly, we consider a distributionally robust chance constrained program (DRCCP) of the following form:
\begin{equation}
\label{eq_drccp_org_ch2}
v^* = \min_{\bm{x} \in \X} \left\{ \bm{c}^\top \bm{x} \colon \inf_{\Pr \in \P} \Pr \left\{ \trxi \colon g(\bm{x}, \trxi) \leq 0 \right\} \geq 1 - \varepsilon \right\},
\end{equation}
which seeks a decision $\bm{x}$ from the feasible region $\X \subseteq \mathbb{R}^n$ that minimizes the cost $\bm{c}^\top \bm{x}$ while ensuring that the constraints $g(\bm{x}, \trxi) \leq 0$ hold with probability at least $1 - \varepsilon$ for all distributions in an ambiguity set $\P$, representing a collection of plausible candidates for the unknown true distribution of $\trxi$. The risk level $\varepsilon \in (0,1)$ specifies the allowable probability of constraint violation, with smaller values of $\varepsilon$ corresponding to more stringent reliability requirements.

To formally characterize $\P$, we employ Wasserstein distance as a metric to quantify the proximity between probability distributions. Specifically, we define the ambiguity set \(\mathcal{P}\) as the set of all probability measures $\Pr$ supported on $\Xi$ whose distance to a reference distribution $\hat{\Pr}$ (e.g., the empirical distribution) is bounded by a radius \(\theta \geq 0\): 
\begin{equation*}
\P = \left\{ \Pr \colon \Pr \{ \trxi \in \Xi \} = 1, W(\Pr, \hat{\Pr}) \leq \theta \right\},
\end{equation*}
where the Wasserstein ambiguity set is defined with type-1 Wasserstein distance: %
\begin{equation*}
W(\Pr_1, \Pr_2) = \inf \left\{ \int_{\Xi \times \Xi} \| \rxi_1 - \rxi_2 \|_p \ \mathbb{Q}(d\rxi_1, d\rxi_2) \colon 
\begin{aligned}
& \mathbb{Q} \text{ is a joint distribution of } \tilde{\rxi}_1 \text{ and } \tilde{\rxi}_2 \\
& \text{with marginals } \Pr_1 \text{ and } \Pr_2, \text{ respectively}
\end{aligned}
\right\},
\end{equation*}
where $\|\cdot\|_p$ with $p\in[1,\infty]$ denotes the $p$-norm on $\Re^m$, so that $\|\rxi_1-\rxi_2\|_p$ represents the cost of transporting a unit of probability mass from $\rxi_1$ to $\rxi_2$; accordingly, $W(\Pr_1,\Pr_2)$ is the minimum expected cost of transporting the probability mass of $\Pr_1$ to that of $\Pr_2$. We denote by $\|\cdot\|_{p^*}$ the dual norm of $\|\cdot\|_p$, i.e., $\|\bm a\|_{p^*}\coloneqq\sup\{\bm a^\top\rxi\colon \|\rxi\|_p\leq 1\}$, where $p^*\in[1,\infty]$ satisfies $1/p+1/p^*=1$ (with the convention $1/\infty=0$).

Following recent studies on Wasserstein DRCCPs \citep[see, e.g.,][]{xie2021distributionally,chen2024data}, we make the following standing assumptions. The deterministic set $\X$ is nonempty, compact, and convex, and the random vector $\trxi$ has support $\Xi=\Re^m$. The constraint function $g$ is the pointwise maximum of $I$ bi-affine functions, that is,
\(
g(\bm{x},\rxi)\coloneqq\max_{i\in[I]} g_i(\bm{x},\rxi)\), with 
\(g_i(\bm{x},\rxi)\coloneqq \bm a_i(\bm x)^\top\rxi-b_i(\bm x)\),
where $\bm a_i\colon\Re^n\to\Re^m$ and $b_i\colon\Re^n\to\Re$ are affine mappings for each $i\in[I]$. When $I=1$, the chance constraint in \eqref{eq_drccp_org_ch2} is an individual chance constraint; when $I>1$, it is a joint chance constraint. We further focus on instances whose coefficient vectors share a common dual norm, namely,
\( \|\bm{a}_i(\bm{x})\|_{p^*}=\|\bm{a}_1(\bm{x})\|_{p^*}\) for all \(\bm{x}\in\X\) and \(i\in[I]\), which holds for individual chance constraints ($I=1$). It also holds when the coefficient vectors do not depend on $\bm x$, i.e., $\bm a_i(\bm x)\equiv\bm a_i$ for all $i\in[I]$ (which includes the case of right-hand-side uncertainty), because each row can then be normalized by the positive constant $1/\|\bm a_i\|_{p^*}$ without changing the feasible region. It holds as well for joint chance constraints with a common coefficient vector, $\bm a_i(\bm x)=\bm a(\bm x)$ for all $i\in[I]$ \citep[see, e.g.,][]{xie2021distributionally}. When the Wasserstein radius $\theta=0$, the Wasserstein ambiguity set reduces to the reference distribution $\hat\Pr$, and Problem~\eqref{eq_drccp_org_ch2} reduces to the regular chance constrained setting. We then assume that \(\theta>0\). Finally, we assume without loss of generality that $b_i(\bm x)\geq 0$ for all $i\in[I]$ whenever $\bm x\in\X$ satisfies $\|\bm{a}_1(\bm{x})\|_{p^*}=0$. Under these assumptions, it follows from Corollary~1 of \citet{xie2021distributionally} (see also Theorem~2 of \citealt{chen2024data}) that DRCCP~\eqref{eq_drccp_org_ch2} admits the following equivalent reformulation:
\begin{align}
\label{1_wasser_reformualtion}
v^*=\min_{\bm x\in\X}\left\{ \bm c^\top \bm x\colon \begin{array}{l}
\displaystyle
\theta\varepsilon^{-1}\left\|\bm{a}_1(\bm{x})\right\|_{p^*}+\hat{\Pr}\text{-}\CVaR_{1-\varepsilon}\left[-\min_{i\in [I]}\left({b}_i(\bm{x})-\bm a_i(\bm x)^\top\hat{\rxi}\right)_+\right]\leq 0
\end{array}
\right\},
\end{align}
where \(\hat{\Pr}\text{-}\CVaR_{1-\varepsilon}(\tilde{X})\) represents the conditional value-at-risk (CVaR) of a random variable \(\tilde{X}\) under distribution \(\hat{\Pr}\) at level \(1-\varepsilon\), as defined by \cite{rockafellar2002conditional}.

\subsection{Literature Review}\label{sec_lit}

Research on DRCCPs has grown rapidly in recent years \citep[see][ for a recent review]{kuccukyavuz2022chance}, wherein the selection of an ambiguity set dictates both the statistical guarantees and the computational complexity of the resulting models. Early studies focused primarily on moment-based ambiguity sets, which encompass all distributions consistent with prescribed moments such as the mean and covariance, because they often admit tractable conic reformulations \citep{zymler2013distributionally,hanasusanto2017ambiguous,xie2017distributionally}. However, such ambiguity sets also contain extreme worst-case distributions, and optimizing against the worst of them can lead to excessive conservatism \citep{gao2023distributionally}. Furthermore, such conservatism persists even as sample sizes increase, as additional data only improve the moment estimates without removing these extreme distributions from the ambiguity set.

Wasserstein ambiguity sets, which have become prevalent in distributionally robust optimization \citep{mohajerin2018data,gao2023distributionally,blanchet2019quantifying,gao2024wasserstein}, circumvent such limitations. Specifically, under mild conditions, the radius of a Wasserstein ambiguity set centered at the empirical distribution can be chosen so that the ambiguity set contains the true distribution with a prescribed confidence level, and this calibrated radius diminishes to zero as the sample size grows \citep{mohajerin2018data}. Consequently, every feasible solution to the DRCCP satisfies the chance constraint under the true distribution with at least a prescribed confidence level, and the model becomes progressively less conservative as more data become available. Due to these benefits, Wasserstein-based methods have become increasingly popular in recent data-driven DRCCP studies \citep{xie2020bicriteria,xie2021distributionally,chen2024data,ji2021data,ho2022distributionally,ho2023strong,shen2023chance,shen2025convex,chen2023approximations}.

Despite these statistical advantages, Wasserstein DRCCPs are computationally difficult to solve \citep{xie2021distributionally,ji2021data,chen2024data,shen2025convex,chen2023approximations}. This difficulty stems in part from the chance constraint itself, because a Wasserstein DRCCP contains a regular chance constrained program as a special case. When the radius is zero, the Wasserstein ambiguity set contains only the empirical distribution, and a Wasserstein DRCCP reduces to a regular chance constrained program under a finite discrete distribution. Any method that solves Wasserstein DRCCPs in general must therefore also solve these chance constrained programs, which are already nonconvex and NP-hard~\citep{luedtke2010integer} even when uncertain parameters only appear on the right-hand side of the constraints. In this sense, Wasserstein DRCCPs inherit the nonconvexity and NP-hardness of regular chance constrained programs, and \citet{xie2020bicriteria} further proved that they remain NP-hard for covering problems.

To address such computational difficulty, existing approaches mainly fall into two categories, exact mixed-integer programming (MIP) reformulations and convex approximations~\citep{kuccukyavuz2022chance}. When the reference distribution is the empirical distribution, a Wasserstein DRCCP admits an exact MIP reformulation \citep{xie2021distributionally,chen2024data}, but the number of binary variables in this reformulation grows with the sample size, which limits its scalability even with strengthened formulations \citep{ho2022distributionally,ho2023strong} and variable fixing \citep{jiang2025terminator}. Convex approximations, on the other hand, sacrifice optimality guarantees for scalability by avoiding binary variables altogether. The CVaR approximation~\citep{nemirovski2007convex} replaces the chance constraint with a convex CVaR constraint, yielding a tractable problem whose solutions are always feasible for the DRCCP but often overly conservative \citep{chen2023approximations}. To alleviate such conservatism, ALSO-X\#~\citep{jiang2025also} performs bisection on the objective value and solves a CVaR-type convex program at each step to find a less conservative solution that remains feasible for the DRCCP. Its computational cost therefore grows with the number of bisection steps. However, ALSO-X\# provides no bound on the number of steps needed to improve on the CVaR solution.

Sequential convex approximation (SCA) methods offer a promising middle ground between these two categories. Rather than searching over the objective value, SCA methods reformulate or approximate the chance constraint as a difference-of-convex (DC) constraint and, at each iteration, convexify it by linearizing its concave part at the current solution, so that each iteration solves a convex program comparable in size to a single CVaR approximation. For regular chance constrained programs, the seminal work of \citet{hong2011sequential} showed that the SCA iterates converge to a Karush--Kuhn--Tucker (KKT) point under mild conditions, and \citet{hong2014conditional} showed that SCA returns a solution at least as good as that of the CVaR approximation. Furthermore, recent studies have proposed proximal and penalty-based DC algorithms for chance constrained programs \citep{wang2026proximal,li2026solving}, and SCA methods have also been applied in engineering settings such as stochastic coordinated beamforming under uncertain channel state information \citep{shi2014optimal} and electric vehicle dispatch under uncertain driver responses \citep{zhang2024joint}.

However, existing SCA methods have largely been limited to regular chance constrained programs, because they rely on DC reformulations or approximations built for a single, fixed distribution \citep{hong2011sequential,hong2014conditional,wang2026proximal,li2026solving}, whereas Wasserstein DRCCPs require the chance constraint to hold for every distribution in the Wasserstein ambiguity set. Extending SCA to Wasserstein DRCCPs therefore requires a fundamentally different and exact DC reformulation. Closely related exact reformulations are available \citep{shen2025convex,jiang2025also}, but, to the best of our knowledge, they have not been cast as a single DC constraint and exploited within an SCA scheme. We develop such a reformulation by showing that the constraint function of a Wasserstein DRCCP equals that of its CVaR approximation minus a convex function that captures the conservatism of the CVaR approximation, yielding an exact DC reformulation of the DRCCP. Building on this reformulation, we develop an SCA algorithm with convergence and iteration-complexity guarantees that, under mild conditions, provably improves on the CVaR approximation after only a few iterations.

\subsection{Motivation}
\label{sec_motivation}
Convex approximations, including the CVaR approximation and refinements such as ALSO-X\#, are widely used in practice~\citep[see, e.g.,][]{hota2019data,xie2021distributionally,chen2023approximations,chen2024data,jiang2025also} for Wasserstein DRCCPs, which are nonconvex and NP-hard in general. These approximations are tractable and always return feasible solutions for the DRCCP, but they can be conservative. The following covering example illustrates this conservatism and the reason behind it, which motivates the proposed DC~reformulation.

\begin{example}
\label{example_sensitivity}
Consider a DRCCP~\eqref{eq_drccp_org_ch2} whose reference distribution \(\hat{\Pr}\) is the empirical distribution of $N=5$ scenarios and $p=2$.
Let 
the risk level be $\varepsilon=0.5$, Wasserstein radius $\theta=0.01$, cost vector $\bm{c} = (1,1,1)^\top$, deterministic feasible region $\X=\{\bm x\in\Re^3_+:\sum_{i\in[3]}x_i\geq 0.01, 0\leq x_i\leq 10000, \forall i\in[3]\}$. The constraint functions are given by $\bm a_i(\x)=\x, b_i(\x)=b =-1, I=1$ with scenarios $\hat{\rxi}^1 = (-1,0,0)^\top$, $\hat{\rxi}^2 =(0,-\alpha,0)^\top$, $\hat{\rxi}^3 = (0,0,-1)^\top$, $\hat{\rxi}^4 = \hat{\rxi}^5 = (-1,-1,-1)^\top$, where scenario $\hat{\rxi}^2$ is perturbed through the parameter $\alpha\in[0.7, 0.9]$ and thus defines a family of instances on which we compare different approximation approaches for solving DRCCPs. 
The resulting DRCCP instances share the same optimal value $10/9$ and the same optimal solution $\bm{x}^* = (10/9,0,0)^\top$. 
However, both the CVaR approximation and ALSO-X\# return feasible solutions whose objective values are far from the optimal value across the entire family, as shown in Figure~\ref{fig:sensitivity}. The CVaR approximation yields objective values of about $3.4$ for every $\alpha$, roughly three times the optimal value. ALSO-X\# improves on the CVaR approximation when $\alpha$ is small but remains conservative, and this improvement vanishes as $\alpha$ increases; for $\alpha\geq 0.76$, ALSO-X\# returns the same objective value as the CVaR approximation.
\QEDB
\vspace{1em}
\begin{figure}[htbp]
\centering
\includegraphics[width=0.6\textwidth]{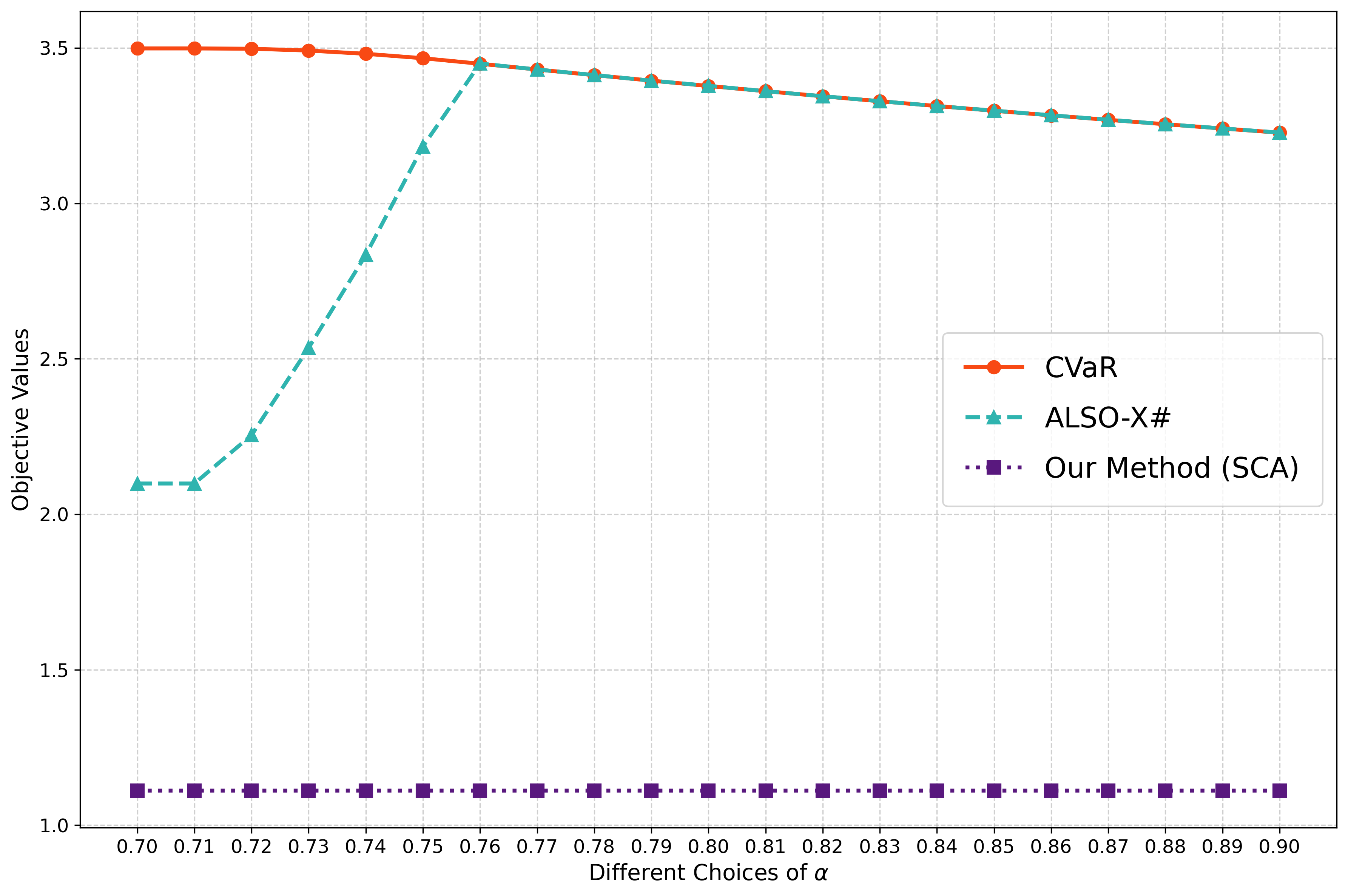}
\caption{Objective values with respect to $\alpha$ returned by the CVaR approximation, ALSO-X\#, and the proposed SCA approach. \textcolor{blue}{}}
\label{fig:sensitivity}
\end{figure}
\end{example}

\begin{figure}[!ht]
\centering
\begin{tikzpicture}[x=1.25cm,y=6cm]

\def\mu{0.5}
\def\sig{1.1}
\def\xmin{-3.2}
\def\xmax{3.5}
\def\varx{-0.35}   
\def\zerox{0.95}   

\def\pdf(#1){ (1/(\sig*sqrt(2*pi))) * exp(-((#1-\mu)^2)/(2*\sig*\sig)) }

\draw[very thick,->] (\xmin,0) -- (\xmax+0.2,0);

\fill[gray!45]
plot[domain=\zerox:\xmax, samples=150] (\x,{\pdf(\x)})
-- (\xmax,0) -- (\zerox,0) -- cycle;

\fill[pattern=north east lines, pattern color=black!65]
plot[domain=\varx:\zerox, samples=150] (\x,{\pdf(\x)})
-- (\zerox,0) -- (\varx,0) -- cycle;

\draw[very thick]
plot[domain=\xmin:\xmax, samples=200] (\x,{\pdf(\x)});

\draw[densely dashed, ultra thick] (\zerox,0) -- (\zerox,0.52);
\node[below] at (\zerox,0) {$0$};

\draw[ultra thick] (\varx,0) -- (\varx,0.52);

\node at (0.3,0.15)  {\small Region II};
\node at (2.0,0.20)  {\small Region III};

\draw[decorate,decoration={brace,amplitude=6pt}]
($( \varx,0.46)$) -- ($( \xmax,0.46)$)
node[midway, yshift=9pt] {\small Region I };

\end{tikzpicture}
\caption{An intuitive decomposition of the tail. Region II represents the part associated with the original Wasserstein DRCCP, while Region III represents the additional positive tail part included by the CVaR approximation.}
\label{fig:tail_decomposition}
\end{figure}
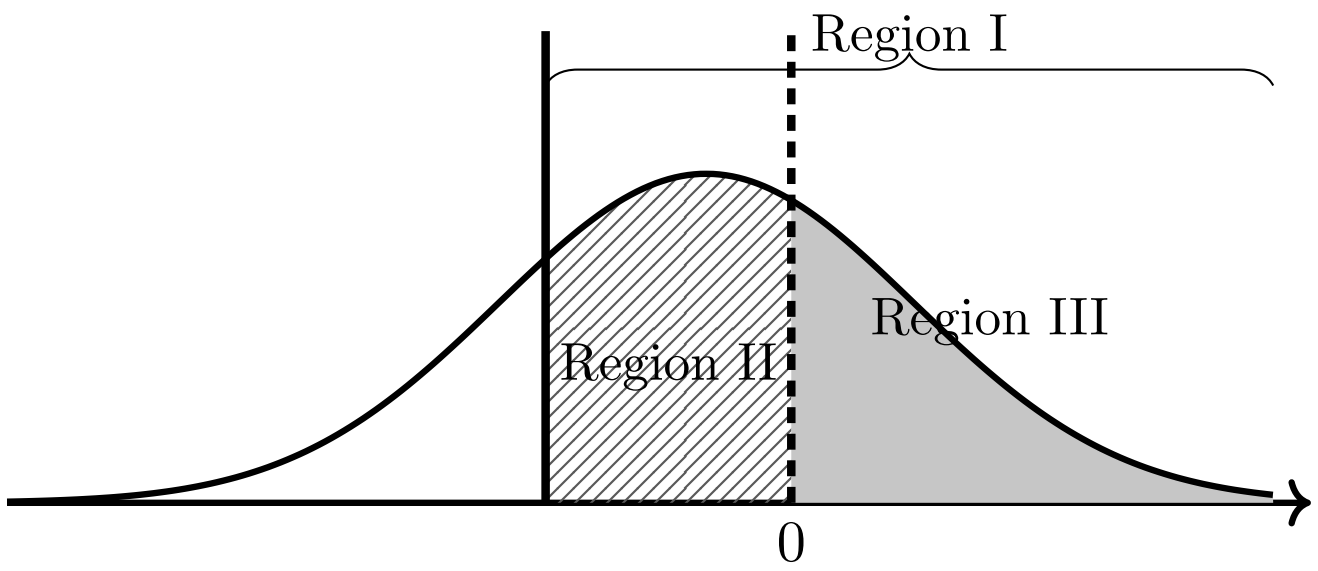

Figure~\ref{fig:tail_decomposition} explains the source of this conservatism by depicting, for a fixed decision $\bm x$, the distribution of the constraint value $g(\bm x,\hat{\rxi})$ under $\hat{\Pr}$. The CVaR approximation restricts the average of $g(\bm x,\hat{\rxi})$ over Region~I, which contains the values above its $(1-\varepsilon)$-quantile and splits at zero into Region~II, whose nonpositive values measure the slack of the constraint in scenarios where it holds, and Region~III, whose positive values measure the size of the violations. The exact reformulation~\eqref{1_wasser_reformualtion} applies CVaR to the truncated value $\min\{g(\bm x,\hat{\rxi}),0\}$, which assigns zero to every scenario in Region~III regardless of how large its violation is, whereas the CVaR approximation applies CVaR to $g(\bm x,\hat{\rxi})$ itself and thus retains the size of each violation. The violations in Region~III therefore account for the conservatism of the CVaR approximation relative to the DRCCP.
This explains the behavior observed in Example~\ref{example_sensitivity}. At the optimal solution $\bm x^*=(10/9,0,0)^\top$, the constraint is violated under scenario $\hat{\rxi}^2$, since $g(\bm x^*,\hat{\rxi}^2)=1>0$ for every $\alpha$, yet $\bm x^*$ remains feasible for the DRCCP, whose constraint evaluates this scenario through the truncated value $\min\{g(\bm x^*,\hat{\rxi}^2),0\}=0$ and thus disregards the size of the violation. The CVaR approximation, in contrast, retains the full size of this violation in Region~III, and excludes $\bm x^*$ from its feasible region. ALSO-X\# is affected in the same way through its CVaR-based lower-level problem. As a result, neither method attains the optimal value $10/9$, and the objective value of ALSO-X\# further varies with $\alpha$ (see the details in Figure~\ref{fig:sensitivity}).

Example~\ref{example_sensitivity} and Figure~\ref{fig:tail_decomposition} suggest that reducing conservatism requires removing the effect of Region~III while retaining the tractability of convex approximations. This observation motivates the DC reformulation developed in this paper, which expresses the constraint of DRCCP~\eqref{1_wasser_reformualtion} as the CVaR approximation minus a convex function that captures Region~III (Theorem~\ref{theorem_DC_equ}). The proposed SCA algorithm replaces this convex function with its linearization at the current solution, which yields a convex inner approximation of the DRCCP of the same size as the CVaR approximation, and updates the linearization at each iteration. As shown by the dotted line with square markers in Figure~\ref{fig:sensitivity}, the SCA algorithm attains the optimal value $10/9$ for every $\alpha\in[0.7,0.9]$. A more detailed analysis is deferred to Section~\ref{sec_alsox}, where we revisit Example~\ref{example_sensitivity} (see also Example~\ref{example_sensitivity_detail} in Appendix~\ref{appendix_example}).

\subsection{Summary of Contributions}
Motivated by these observations, we now summarize the main contributions of this paper:
\begin{itemize}
  \item[(i)] We derive an equivalent difference-of-convex (DC) reformulation of Wasserstein DRCCPs, which decomposes the distributionally robust chance constraint into the difference of two convex functions. It shows that the CVaR approximation drops a nonnegative convex term, which explains its conservatism and suggests how to construct tighter convex inner approximations;

\item[(ii)] Exploiting this DC structure, we develop a sequential convex approximation (SCA) algorithm whose iterates remain feasible for the DRCCP and whose objective values improve monotonically. We show that these iterates converge to critical points of the DC reformulation and derive an iteration-complexity bound for reaching an approximate critical point;

\item[(iii)] We further compare the SCA algorithm with the CVaR approximation and ALSO-X\#~\citep{jiang2025also}, the state-of-the-art convex approximation for Wasserstein DRCCPs. We prove that, under a simple initialization rule, the SCA algorithm returns a solution no worse than the CVaR approximation after a single iteration. For ALSO-X\#, we construct two instances in which each method attains a strictly better objective value than the other, showing that neither dominates the other;

\item[(iv)] %
We conduct numerical experiments on portfolio optimization problems with individual and joint chance constraints. Across all tested instances, SCA yields lower objective values on average than the CVaR approximation and ALSO-X\#. For individual chance constraints, its objective values match those of the exact Big-M formulation to the reported precision. For joint chance constraints, SCA reduces the average objective value by up to $39.58\%$ relative to CVaR and $27.13\%$ relative to ALSO-X\#. It requires $3.33$--$6.05$ seconds on average and matches or improves the best solutions found by Big-M within a $3600$-second time limit.

\end{itemize}

\noindent\textbf{Organization.} The remainder of this paper is organized as follows. Section~\ref{sec_SCA_algorithm} derives an equivalent DC reformulation of DRCCP~\eqref{1_wasser_reformualtion}, develops the SCA algorithm, and establishes its convergence and iteration complexity. Section~\ref{sec_comparisons} compares the SCA algorithm with the CVaR approximation and ALSO-X\#. Section~\ref{sec_Numerical_study} then reports numerical results for these methods and for an exact mixed-integer benchmark on portfolio problems with individual and joint chance constraints. %

\noindent\textbf{Notation.} We use bold letters (e.g., $\vect{x},\vect{A}$) to denote vectors and matrices and use corresponding non-bold letters to denote their components. We denote realizations of the random vector $\trxi$ by $\rxi$. Given an integer $n$, we let $[n]:=\{1,2,\ldots,n\}$, and use $\Re_+^n:=\{\bm {x}\in \Re^n:x_i\geq0, \forall i\in [n]\}$. Given a real number $t$, we let $(t)_+:=\max\{t,0\}$ and $(t)_-:=\min\{t,0\}$. Given a probability distribution $\Pr$ on $\Xi$, we use $\Pr\{A\}$ to denote $\Pr\{\bm\xi:\text{condition} \ A(\bm\xi) \ \text{holds}\}$ when $A(\bm\xi)$ is a condition on $\bm\xi$, and to denote $\Pr\{\bm\xi\colon \bm\xi \in A\}$ when $A \subseteq \Xi$ is $\Pr$\text{-}measurable.  
For a set $R$, the indicator function $\I(\bm x\in R)$ equals $1$ if $\bm x\in R$ and $0$ otherwise. The normal cone at point $\x\in R$ is denoted by $\N_R(\x):=\{\bm{h}: \bm{h}^\top(\hat{\x}-\x)\leq 0, \forall \hat{\x}\in R\}$, and $\emptyset$ if $\x\notin R$. For a convex function $f\colon\Re^n\to\Re$, $\partial f(\x):=\{\bm s\in\Re^n\colon f(\bm y)\geq f(\x)+\bm s^\top(\bm y-\x),\ \forall \bm y\in\Re^n\}$ denotes its subdifferential at $\x$, and for two sets $A,B\subseteq\Re^n$, $A-B:=\{\bm a-\bm b\colon \bm a\in A,\ \bm b\in B\}$. We use the phrase ``better than" to indicate ``at least as good as." %

\section{An Equivalent DC Reformulation of DRCCP and the SCA Algorithm}
\label{sec_SCA_algorithm}

In this section, we formalize the tail decomposition illustrated in Figure~\ref{fig:tail_decomposition}, use it to derive an equivalent DC reformulation of DRCCP~\eqref{1_wasser_reformualtion}, and then develop an SCA algorithm based on this reformulation. We begin by introducing the risk measures that underlie the analysis. For a scalar random variable $\tilde X$ with distribution $\hat{\Pr}$ and cumulative distribution function $F_{\tilde X}(s):=\hat{\Pr}\{\tilde X\le s\}$, the $u$-value-at-risk (VaR) of $\tilde X$ is its $u$-quantile $\hat{\Pr}\text{-}\VaR_u(\tilde X):=\inf\{s\in\Re\colon F_{\tilde X}(s)\geq u\}$ for $u\in(0,1)$, and the $(1-\varepsilon)$-conditional value-at-risk (CVaR) of $\tilde X$ is
\begin{equation}
\label{eq_cvar_def}
\hat{\Pr}\text{-}\CVaR_{1-\varepsilon}(\tilde{X}) := \min_{\beta\in\Re} \left\{ \beta + \frac{1}{\varepsilon}\, \E_{\hat\Pr} \big[(\tilde{X} - \beta)_+\big] \right\}.
\end{equation}
Furthermore, by \citet[ Proposition~3.2 and Corollary~4.3]{acerbi2002coherence}, CVaR admits an equivalent quantile representation 
\begin{equation}
\label{eq_cvar_quantile}
\hat{\Pr}\text{-}\CVaR_{1-\varepsilon}(\tilde{X}) = \frac{1}{\varepsilon}\int_{1-\varepsilon}^{1}\hat{\Pr}\text{-}\VaR_u(\tilde X)\,du,
\end{equation}
which shows that the $(1-\varepsilon)$-CVaR of $\tilde X$ can be viewed as the expected value of $\tilde X$ over its upper $\varepsilon$-tail $\{\tilde X\geq\hat{\Pr}\text{-}\VaR_{1-\varepsilon}(\tilde X)\}$. 

This quantile representation provides a transport interpretation of the deterministic reformulation~\eqref{1_wasser_reformualtion} \citep{chen2024data}. We call a scenario $\rxi$ unsafe for decision $\bm x$ if $\bm x$ violates the constraint under $\rxi$, that is, $g(\bm x,\rxi)>0$. The infimum in DRCCP~\eqref{eq_drccp_org_ch2} can then be viewed as an adversary that transports probability mass of $\hat{\Pr}$ onto unsafe scenarios, subject to a total transportation cost of at most $\theta$. To make a probability mass of $\varepsilon$ unsafe as cheaply as possible, the adversary leaves the scenarios that are already unsafe at $\bm x$ in place at zero cost and moves only the currently safe scenarios that are closest to the unsafe set, each at a cost equal to its distance to that set, proportional to $\min_{i\in[I]}\big(b_i(\bm x)-\bm a_i(\bm x)^\top\hat{\rxi}\big)_+=-g(\bm x,\hat{\rxi})_-$. The deterministic reformulation~\eqref{1_wasser_reformualtion} then shows that the minimum cost of making a probability mass of $\varepsilon$ unsafe is the conditional expected distance to the unsafe set over the $\varepsilon$-fraction of scenarios closest to it, and a decision $\bm x$ is feasible exactly when the budget $\theta$ is insufficient to make more than a probability mass of $\varepsilon$ unsafe.

Although treating every unsafe scenario as costless is faithful to the chance constraint, which limits only the probability of violation and not its magnitude, it also makes DRCCP~\eqref{1_wasser_reformualtion} nonconvex in general, as $g(\bm x,\hat{\rxi})$ is convex in $\bm x$ whereas $g(\bm x,\hat{\rxi})_-$ is not. To address this, the CVaR approximation \citep{nemirovski2007convex} replaces $g(\bm x,\hat{\rxi})_-$ with $g(\bm x,\hat{\rxi})$, yielding the following convex program
\begin{align}
v^{\CVaR}=\min_{\bm x\in\X}\left\{ \bm c^\top \bm x\colon\theta\varepsilon^{-1}\left\|\bm{a}_1(\bm{x})\right\|_{p^*} + \hat{\Pr}\text{-}\CVaR_{1-\varepsilon}\Big[ \max_{i\in[I]} \big\{\bm{a}_i(\bm{x})^\top\hat{\rxi}-b_i(\bm{x})\big\}\Big]\leq 0 \right\},\label{worse_case_cvar_q_formulation_a}
\end{align}
where unsafe scenarios now enter with a negative distance $-g(\bm x,\hat{\rxi})/\|\bm a_1(\bm x)\|_{p^*}$, as if the adversary were credited in proportion to how severely the scenario violates the constraint. Since $g(\bm x,\hat{\rxi})_-\leq g(\bm x,\hat{\rxi})$, this credit can only help the adversary, so every feasible solution of~\eqref{worse_case_cvar_q_formulation_a} is feasible for DRCCP~\eqref{1_wasser_reformualtion} and $v^{\CVaR}\geq v^*$. The decomposition in Figure~\ref{fig:tail_decomposition} visually illustrates the relationship between the original DRCCP and its CVaR approximation. Specifically, Region~I represents the upper $\varepsilon$-tail of $g(\bm x,\hat{\rxi})$ and corresponds to the $\hat{\Pr}\text{-}\CVaR_{1-\varepsilon}[g(\bm{x},\hat{\rxi})]$ in~\eqref{worse_case_cvar_q_formulation_a}. The safe scenarios in this tail, namely Region~II, correspond to the $\hat{\Pr}\text{-}\CVaR_{1-\varepsilon}[g(\bm{x},\hat{\rxi})_-]$ in DRCCP~\eqref{1_wasser_reformualtion}, capturing the adversary's minimum transportation cost. The unsafe scenarios in this tail, namely Region~III, correspond to $\hat{\Pr}\text{-}\CVaR_{1-\varepsilon}[g(\bm{x},\hat{\rxi})_+]$, the credit that the CVaR approximation grants the adversary and also the source of the over-conservatism.

To reduce this conservatism while retaining the tractability of the CVaR approximation, we establish in Section~\ref{section_dc_reform} the additive decomposition 
\(\hat{\Pr}\text{-}\CVaR_{1-\varepsilon}\big[g(\bm x,\hat{\rxi})\big] = \hat{\Pr}\text{-}\CVaR_{1-\varepsilon}\big[g(\bm x,\hat{\rxi})_-\big] + \hat{\Pr}\text{-}\CVaR_{1-\varepsilon}\big[g(\bm x,\hat{\rxi})_+\big]\), which motivates the equivalent DC reformulation
\begin{align}
\label{eq_dc_orginial}
v^*(\mathrm{DC})=\min_{\bm x\in\X}\left\{ \bm c^\top \bm x\colon 
\begin{aligned}
&\theta\varepsilon^{-1}\left\|\bm{a}_1(\bm{x})\right\|_{p^*} + \hat{\Pr}\text{-}\CVaR_{1-\varepsilon}\Big[ \max_{i\in[I]} \big\{\bm{a}_i(\bm{x})^\top\hat{\rxi}-b_i(\bm{x})\big\}\Big] \\
&\quad\quad\leq \hat{\Pr}\text{-}\CVaR_{1-\varepsilon}\Big[ \max_{i\in[I]} \big(\bm{a}_i(\bm{x})^\top\hat{\rxi}-b_i(\bm{x})\big)_+\Big]
\end{aligned}
\right\}.
\end{align}
Subsequently, we propose an SCA algorithm for solving it in Section~\ref{section_sub_sca}.

\subsection{An Equivalent DC Reformulation of DRCCP}
\label{section_dc_reform}
To show the equivalence between DC reformulation~\eqref{eq_dc_orginial} and DRCCP~\eqref{1_wasser_reformualtion}, we first recall the definition of comonotonicity (see, e.g., Proposition~4.5 in \citealt{denneberg1994non}): scalar random variables $\tilde{X}$ and $\tilde{Y}$ defined on a sample space $\Omega$ are comonotonic if there exists a set $\Omega_0\subseteq\Omega$ of probability one such that $\big(\tilde{X}(\omega_1)-\tilde{X}(\omega_2)\big)\big(\tilde{Y}(\omega_1)-\tilde{Y}(\omega_2)\big)\geq 0$ for all $\omega_1,\omega_2\in\Omega_0$.
\begin{lemma}
\label{lemma_comonotonic} (\citealt{denneberg1994non,dhaene2002concept})
Let $\tilde{X}$ and $\tilde{Y}$ be two comonotonic random variables. Then for all $\varepsilon \in (0,1)$, it holds that
\[
\hat{\Pr}\text{-}\mathrm{VaR}_{1-\varepsilon}(\tilde{X} + \tilde{Y}) = \hat{\Pr}\text{-}\mathrm{VaR}_{1-\varepsilon}(\tilde{X}) + \hat{\Pr}\text{-}\mathrm{VaR}_{1-\varepsilon}(\tilde{Y}).
\]
\end{lemma}
From Lemma~\ref{lemma_comonotonic}, we can show that
$\hat{\Pr}\text{-}\CVaR_{1-\varepsilon}[\tilde{X}] = \hat{\Pr}\text{-}\CVaR_{1-\varepsilon}[\tilde{X}^+] + \hat{\Pr}\text{-}\CVaR_{1-\varepsilon}[\tilde{X}^-]$ for any scalar random variable $\tilde{X}$ and $\varepsilon \in (0,1)$, 
where $\tilde{X}^+ = \max\{\tilde{X}, \bm0\}$ and $\tilde{X}^- = \min\{\tilde{X}, \bm0\}$.
This helps us establish the equivalence between DC reformulation \eqref{eq_dc_orginial} and DRCCP~\eqref{1_wasser_reformualtion}.

\begin{theorem}
\label{theorem_DC_equ}
DC reformulation \eqref{eq_dc_orginial} is equivalent to DRCCP~\eqref{1_wasser_reformualtion}, i.e., $v^*(\mathrm{DC}) = v^*$. 
\end{theorem}
\begin{subequations}
\begin{proof}
We split the proof into three parts. 

\noindent\textbf{Part 1.}
We first show that $\hat{\Pr}\text{-}\CVaR_{1-\varepsilon}[\tilde{X}] = \hat{\Pr}\text{-}\CVaR_{1-\varepsilon}[\tilde{X}^+] + \hat{\Pr}\text{-}\CVaR_{1-\varepsilon}[\tilde{X}^-]$ for a given random variable $\tilde{X}$ and $\varepsilon \in (0,1)$.
Since we know that $\tilde{X}^+$ and $\tilde{X}^-$ are both monotone increasing, according to the results in Lemma~\ref{lemma_comonotonic}, we have:
\begin{align*}
\hat{\Pr}\text{-}\CVaR_{1-\varepsilon}[\tilde{X}] &= \frac{1}{\varepsilon} \int_{1-\varepsilon}^1 \hat{\Pr}\text{-}\mathrm{VaR}_{u}(\tilde{X})du = \frac{1}{\varepsilon} \int_{1-\varepsilon}^1 \hat{\Pr}\text{-}\mathrm{VaR}_{u}(\tilde{X}^+ +\tilde{X}^-)du \\
&= \frac{1}{\varepsilon} \int_{1-\varepsilon}^1 \left[ \hat{\Pr}\text{-}\mathrm{VaR}_{u}(\tilde{X}^+) + \hat{\Pr}\text{-}\mathrm{VaR}_{u}(\tilde{X}^-) \right]du \\
& = \frac{1}{\varepsilon} \int_{1-\varepsilon}^1 \hat{\Pr}\text{-}\mathrm{VaR}_{u}(\tilde{X}^+)du + \frac{1}{\varepsilon} \int_{1-\varepsilon}^1 \hat{\Pr}\text{-}\mathrm{VaR}_{u}(\tilde{X}^-)du \\
&= \hat{\Pr}\text{-}\CVaR_{1-\varepsilon}[\tilde{X}^+] + \hat{\Pr}\text{-}\CVaR_{1-\varepsilon}[\tilde{X}^-].
\end{align*}
Hence, for any random variable $\tilde{X}$ and risk level $\varepsilon \in (0,1)$, we have
\begin{equation}
\hat{\Pr}\text{-}\CVaR_{1-\varepsilon}[\tilde{X}] = \hat{\Pr}\text{-}\CVaR_{1-\varepsilon}[\tilde{X}^+] + \hat{\Pr}\text{-}\CVaR_{1-\varepsilon}[\tilde{X}^-].
\label{CVaR_comonotonicity}
\end{equation}

\noindent{\textbf{Part 2.}}
Given finitely many random variables $\{\tilde{X}_i\}$, define $\tilde{U} = \max_i \{\tilde{X}_i\}$, $\tilde{U}^+ = \max_i \{\tilde{X}_i^+\}$, and $\tilde{U}^- = \max_i \{\tilde{X}_i^-\}$, where $\tilde{X}_i^+ = \max\{\tilde{X}_i, 0\}$ and $\tilde{X}_i^- = \min\{\tilde{X}_i, 0\}$.  If $\tilde{U}(\omega)\geq 0$ for some $\omega\in\Omega$, then $\tilde{U}^+(\omega)=\tilde{U}(\omega)$ and $\tilde{U}^-(\omega)=0$. If $\tilde{U}(\omega)< 0$, then $\tilde{U}^+(\omega)=0$ and $\tilde{U}^-(\omega)=\tilde{U}(\omega)$. In either case, we have $\tilde{U}(\omega)=\tilde{U}^+(\omega)+\tilde{U}^-(\omega)$, establishing that $\tilde{U} = \tilde{U}^+ + \tilde{U}^-$ almost surely. Applying the equality \eqref{CVaR_comonotonicity}, we obtain:
\begin{equation}
\hat{\Pr}\text{-}\CVaR_{1-\varepsilon}\left[\max\limits_i \tilde{X}_i\right] = \hat{\Pr}\text{-}\CVaR_{1-\varepsilon}\left[\max\limits_i \{\tilde{X}_i^+\}\right] + \hat{\Pr}\text{-}\CVaR_{1-\varepsilon}\left[\max\limits_i \{\tilde{X}_i^-\}\right].
\label{max_CVaR_comonotonicity}
\end{equation}

\noindent{\textbf{Part 3.}}
Substituting $\bm{a}_i(\bm{x})^\top \widehat{\bm{\xi}} - b_i(\bm{x})= \tilde{X}_i$ demonstrates the equivalence between DC reformulation \eqref{eq_dc_orginial} and DRCCP~\eqref{1_wasser_reformualtion}. This completes the proof.
\QEDA
\end{proof}
\end{subequations}
We make the following remarks on Theorem~\ref{theorem_DC_equ}:
\begin{itemize}
\item[(i)] The equivalence holds for any reference distribution $\hat{\Pr}$ with a finite first moment; and
\item[(ii)] The formulation of DRCCP~\eqref{eq_dc_orginial} can be used to alleviate the over-conservatism of CVaR approximation~\eqref{worse_case_cvar_q_formulation_a}, as formally discussed in the next section.
\end{itemize}
According to Theorem~\ref{theorem_DC_equ}, we know that solving the DC reformulation~\eqref{eq_dc_orginial} is, in general, NP-hard (see, e.g., \cite{xie2020bicriteria}). Rather than pursuing exact methods, we focus on efficient approaches that provide high-quality upper bounds for DRCCP~\eqref{eq_dc_orginial}. To begin with, the right-hand side of the constraint in DRCCP~\eqref{eq_dc_orginial} can be simplified further. The following result adapts Corollary~1 of \citet{shen2025convex} to our setting.
\begin{theorem}
\label{theorem_DC_equ_simple}
DRCCP~\eqref{eq_dc_orginial} is equivalent to   
\begin{align}
\label{eq_dc_orginial_simple}
v^*(\mathrm{DC})=\min_{\bm x\in\X}\left\{ \bm c^\top \bm x\colon 
\begin{aligned}
&\theta\varepsilon^{-1}\left\|\bm{a}_1(\bm{x})\right\|_{p^*} + \hat{\Pr}\text{-}\CVaR_{1-\varepsilon}\left[ \max_{i\in[I]} \left\{\bm{a}_i(\bm{x})^\top\hat{\rxi}-b_i(\bm{x})\right\}\right] \\
&\quad\quad\leq \frac{1}{\varepsilon}\E_{\hat\Pr}\left[ \max_{i\in[I]} \left\{\bm{a}_i(\bm{x})^\top\hat{\rxi}-b_i(\bm{x})\right\}_+\right]
,\\
& \hat{\Pr}\left\{\bm a_i(\bm x)^\top\hat{\rxi}-b_i(\bm{x})\leq 0, \forall i\in[I] \right\}\geq 1-\varepsilon
\end{aligned}
\right\}.
\end{align}
\end{theorem}
\begin{proof}
Based on the definition of CVaR, for any fixed $\bm x\in \X$, we have
\begin{align*}
&\hat{\Pr}\text{-}\CVaR_{1-\varepsilon}
\left[
\left(
\max_{i\in[I]}
\left\{
\bm a_i(\bm x)^\top\hat{\rxi}-b_i(\bm x)
\right\}
\right)_+
\right] \\
&\qquad = \min_{\beta\in\mathbb R}
\left\{
\beta+
\frac{1}{\varepsilon}
\E_{\hat\Pr}
\left[
\left(
\left(
\max_{i\in[I]}
\left\{
\bm a_i(\bm x)^\top\hat{\rxi}-b_i(\bm x)
\right\}
\right)_+
-\beta
\right)_+
\right]
\right\}.
\end{align*}
We first show that, whenever $\hat{\Pr}
\{
\bm a_i(\bm x)^\top\hat{\rxi}-b_i(\bm x)\leq 0,
\ \forall i\in[I]
\}
\geq 1-\varepsilon$, we have
\begin{align*}
\hat{\Pr}\text{-}\CVaR_{1-\varepsilon}
\left[
\left(
\max_{i\in[I]}
\left\{
\bm a_i(\bm x)^\top\hat{\rxi}-b_i(\bm x)
\right\}
\right)_+
\right]
=
\frac{1}{\varepsilon}
\E_{\hat\Pr}
\left[
\left(
\max_{i\in[I]}
\left\{
\bm a_i(\bm x)^\top\hat{\rxi}-b_i(\bm x)
\right\}
\right)_+
\right].
\end{align*}
Indeed, the above chance constraint implies $\hat{\Pr}
\{
(
\max_{i\in[I]}
\{
\bm a_i(\bm x)^\top\hat{\rxi}-b_i(\bm x)
\}
)_+
=0
\}
\geq 1-\varepsilon$. Since we also know 
$
(
\max_{i\in[I]}
\{
\bm a_i(\bm x)^\top\hat{\rxi}-b_i(\bm x)
\}
)_+
\geq 0$,  we obtain
\begin{align*}
\hat{\Pr}\text{-}\mathrm{VaR}_{1-\varepsilon}
\left[
\left(
\max_{i\in[I]}
\left\{
\bm a_i(\bm x)^\top\hat{\rxi}-b_i(\bm x)
\right\}
\right)_+
\right]
=0.
\end{align*}
According to the discussions in Section 3 of
\citet{rockafellar2000optimization}, we know $\beta=0$ is an optimal solution of $\hat{\Pr}\text{-}\CVaR_{1-\varepsilon}
[
(
\max_{i\in[I]}
\{
\bm a_i(\bm x)^\top\hat{\rxi}-b_i(\bm x)
\}
)_+
]$, and hence
\begin{align*}
\hat{\Pr}\text{-}\CVaR_{1-\varepsilon}
\left[
\left(
\max_{i\in[I]}
\left\{
\bm a_i(\bm x)^\top\hat{\rxi}-b_i(\bm x)
\right\}
\right)_+
\right] =
\frac{1}{\varepsilon}
\E_{\hat\Pr}
\left[
\left(
\max_{i\in[I]}
\left\{
\bm a_i(\bm x)^\top\hat{\rxi}-b_i(\bm x)
\right\}
\right)_+
\right].
\end{align*}
Now suppose that $\bm x$ is feasible to
DRCCP~\eqref{eq_dc_orginial}. Since the reference distribution $\hat{\Pr}$
belongs to the Wasserstein ambiguity set, the corresponding regular
chance constraint is satisfied, i.e.,
\begin{align*}
\hat{\Pr}
\left\{
\bm a_i(\bm x)^\top\hat{\rxi}-b_i(\bm x)\leq 0,
\ \forall i\in[I]
\right\}
\geq 1-\varepsilon.
\end{align*}
Thus, the equality above holds. Substituting this equality into
DRCCP~\eqref{eq_dc_orginial} gives the first constraint in
DRCCP~\eqref{eq_dc_orginial_simple}, while the second constraint in
DRCCP~\eqref{eq_dc_orginial_simple} also holds.

Conversely, suppose that $\bm x$ is feasible to
DRCCP~\eqref{eq_dc_orginial_simple}. Its second constraint implies
\begin{align*}
\hat{\Pr}
\left\{
\bm a_i(\bm x)^\top\hat{\rxi}-b_i(\bm x)\leq 0,
\ \forall i\in[I]
\right\}
\geq 1-\varepsilon.
\end{align*}
Therefore, the same equality holds, that is,
\begin{align*}
\hat{\Pr}\text{-}\CVaR_{1-\varepsilon}
\left[
\left(
\max_{i\in[I]}
\left\{
\bm a_i(\bm x)^\top\hat{\rxi}-b_i(\bm x)
\right\}
\right)_+
\right] =
\frac{1}{\varepsilon}
\E_{\hat\Pr}
\left[
\left(
\max_{i\in[I]}
\left\{
\bm a_i(\bm x)^\top\hat{\rxi}-b_i(\bm x)
\right\}
\right)_+
\right].
\end{align*}
Substituting this equality into the first constraint of DRCCP~\eqref{eq_dc_orginial_simple} gives exactly the constraint in DRCCP~\eqref{eq_dc_orginial}. Hence, DRCCP~\eqref{eq_dc_orginial} and DRCCP~\eqref{eq_dc_orginial_simple} are equivalent.
\QEDA
\end{proof}

We make the following remarks on Theorem~\ref{theorem_DC_equ_simple}:
\begin{itemize}
\item[(i)] Following recent studies in chance constrained programs (see, e.g., \citealt{jiang2022also}), the term $\E_{\hat\Pr}[ \max_{i\in[I]} \{\bm{a}_i(\bm{x})^\top\hat{\rxi}-b_i(\bm{x})\}_+]$ can be referred to as the hinge-loss function. If we replace it with $0$ on the right-hand side of the constraint in DRCCP~\eqref{eq_dc_orginial_simple}, we recover CVaR approximation \eqref{worse_case_cvar_q_formulation_a}; 
\item[(ii)] Notably, DRCCP~\eqref{eq_dc_orginial_simple} remains valid for any deterministic set $\X$, e.g., when the decision variables are binary (i.e., $\X = \{0,1\}^n$); and 
\item[(iii)] Both DRCCP~\eqref{eq_dc_orginial} and DRCCP~\eqref{eq_dc_orginial_simple} admit DC representations involving CVaR-type terms. It is worth noting that the difference between CVaR-type terms has appeared in the stochastic programming literature, such as the $\text{In-CVaR}_{\alpha}^{1-\beta}(\cdot)$ notation introduced by \citet{liu2023risk}, which is defined as
\begin{equation*}
  \text{In-CVaR}^{1-\beta}_{\alpha} (g(\bm{x}, \trxi)) := \frac{1-\alpha}{1-\alpha - \beta}\hat{\Pr}\text{-}\CVaR_\alpha (g(\bm{x}, \trxi)) - \frac{\beta}{1-\alpha - \beta}\hat{\Pr}\text{-}\CVaR_{1-\beta} (g(\bm{x}, \trxi)).
\end{equation*}
We refer interested readers to \citet{liu2023risk} for further discussion.
\end{itemize}
Even though Theorem~\ref{theorem_DC_equ_simple} gives an alternative exact formulation,  the constraint $\hat{\Pr}\{\bm a_i(\bm x)^\top\hat{\rxi}-b_i(\bm{x})\leq 0, \forall i\in[I] \}\geq 1-\varepsilon$ may complicate the subsequent algorithmic analysis. Therefore, for ease of analysis, we work with DRCCP~\eqref{eq_dc_orginial}. For notational convenience, we define
\begin{align*}
& g_1(\bm{x}) :=\theta\varepsilon^{-1}\left\|\bm{a}_1(\bm{x})\right\|_{p^*}
+ \hat{\Pr}\text{-}\CVaR_{1-\varepsilon}\!\left[ \max_{i\in[I]}\{\bm{a}_i(\bm{x})^\top\hat{\rxi}-b_i(\bm{x})\}\right],\\
&g_2(\bm{x}) := \hat{\Pr}\text{-}\CVaR_{1-\varepsilon}\left[ \max_{i\in[I]} \left\{\bm{a}_i(\bm{x})^\top\hat{\rxi}-b_i(\bm{x})\right\}_+\right].
\end{align*}
Let $\Omega(\mathrm{DC}):=\{\bm x\in\X\colon g_1(\bm x)\leq g_2(\bm x)\}$ denote the feasible region of DRCCP~\eqref{eq_dc_orginial}. For any $\bm y\in\X$ and $\bm s_2(\bm y)\in\partial g_2(\bm y)$, let
\(
\Omega(\mathrm{DC})_{\bm y}:=\big\{\bm x\in\X\colon g_1(\bm x)\leq g_2(\bm y)+\bm s_2(\bm y)^\top(\bm x-\bm y)\big\}
\)
denote the convex set obtained from $\Omega(\mathrm{DC})$ by replacing $g_2$ with its linearization at $\bm y$. Since $g_1$ and $g_2$ are convex and the objective $\bm c^\top\bm x$ is linear, DRCCP~\eqref{eq_dc_orginial} is a DC program.

\subsection{SCA Algorithm for DRCCP}
\label{section_sub_sca}
In this subsection, we propose an SCA algorithm to solve DRCCP~\eqref{eq_dc_orginial}. The main idea is to exploit the DC structure of its constraint $g_1(\bm x)\leq g_2(\bm x)$ and iteratively solve a sequence of convex approximations. Specifically, at each iteration $k$ of Algorithm~\ref{alg_DCA}, we construct and solve the convex program~\eqref{eq_dc}, which is obtained from DRCCP~\eqref{eq_dc_orginial} by replacing the convex function $g_2$ on the right-hand side of the constraint with its linearization $g_2(\bm x^{(k)})+(\bm s_2^{(k)})^\top(\bm x-\bm x^{(k)})$ around the current iterate $\bm{x}^{(k)}$, where $\bm s_2^{(k)}\in\partial g_2(\bm x^{(k)})$ denotes an arbitrary subgradient of $g_2$ at $\bm x^{(k)}$. The feasible region of~\eqref{eq_dc} is thus $\Omega(\mathrm{DC})_{\bm x^{(k)}}$ with $\bm s_2(\bm x^{(k)})=\bm s_2^{(k)}$. An optimal solution $\bm{x}^{(k+1)}$ of the convex subproblem~\eqref{eq_dc} then becomes the next iterate.
This process is repeated until the change in objective value between consecutive iterates falls below a predefined tolerance $\tau$. Since each subproblem~\eqref{eq_dc} differs from CVaR approximation~\eqref{worse_case_cvar_q_formulation_a} only by an affine right-hand side, it is a convex program of the same size and can be solved by the same solvers as the CVaR approximation. Moreover, since $g_2$ is convex, its linearization is a global lower bound on $g_2$, so $\Omega(\mathrm{DC})_{\bm x^{(k)}}\subseteq\Omega(\mathrm{DC})$. Starting from a feasible solution $\bm x^{(0)}\in\Omega(\mathrm{DC})$, every iterate, including the final iterate $\bm{x}^{\mathrm{DC}}$, is thus feasible for DRCCP~\eqref{eq_dc_orginial}, and its objective value $v^{\mathrm{DC}}$ provides a valid upper bound on the optimal value $v^*$. Algorithm~\ref{alg_DCA} summarizes this SCA procedure.

\begin{algorithm}[H]
\caption{An SCA Algorithm for   DRCCP~\eqref{eq_dc_orginial}, or equivalently,  DRCCP~\eqref{eq_dc_orginial_simple}}
\label{alg_DCA}
\begin{algorithmic}[1]
\State \textbf{Input:} Given $\bm{x}^{(0)} \in \Omega(\mathrm{DC}) $, tolerance $\tau>0$, and set $k=0$ and $v^{(-1)}(\mathrm{DC})=\bm c^\top \bm{x}^{(0)}$
\State Solve SCA formulation \eqref{eq_dc} to obtain its optimal solution $\bm{x}^{(k+1)}$
\begin{align}
\label{eq_dc}
v^{(k)}(\mathrm{DC})=\min_{\bm x\in\X}\left\{ \bm c^\top \bm x\colon 
\begin{aligned}
&\theta\varepsilon^{-1}\left\|\bm{a}_1(\bm{x})\right\|_{p^*} + \hat{\Pr}\text{-}\CVaR_{1-\varepsilon}\left[ \max_{i\in[I]} \left\{\bm{a}_i(\bm{x})^\top\hat{\rxi}-b_i(\bm{x})\right\}\right]\\
& \quad \leq \hat{\Pr}\text{-}\CVaR_{1-\varepsilon}\left[ \max_{i\in[I]} \left\{\bm{a}_i(\bm{x}^{(k)})^\top\hat{\rxi}-b_i(\bm{x}^{(k)})\right\}_+\right]\\
& \qquad {}+ \left(\bm{s}_2^{(k)}\right)^\top \left[ \bm x-\bm x^{(k)}\right]
\end{aligned}
\right\},
\end{align}
where $\bm{s}_2^{(k)} \in \partial g_2(\bm x^{(k)})$,
and $k$ denotes the $k$-th iteration
\State If $|\bm c^\top\bm{x}^{(k+1)} - \bm c^\top\bm{x}^{(k)}| < \tau$, i.e.,  $|v^{(k)}(\mathrm{DC})-v^{(k-1)}(\mathrm{DC})|< \tau$, go to Step 5
\State Set $k=k+1$ and return to Step 2 
\State \textbf{Output:} A solution $\bm x^{\mathrm{DC}}=\bm x^{(k+1)}$ and its objective value $v^{\mathrm{DC}}=\bm c^\top\bm x^{(k+1)}$ to DRCCP~\eqref{eq_dc_orginial}
\end{algorithmic}
\end{algorithm}
Algorithm~\ref{alg_DCA} applies the SCA used in~\citet{hong2011sequential,hong2014conditional} to DRCCP~\eqref{eq_dc_orginial}. In our setting, when $\theta>0$, each iteration of Algorithm~\ref{alg_DCA} solves the SCA formulation with an additional regularization term $\theta\varepsilon^{-1}\left\|\bm{a}_1(\bm{x})\right\|_{p^*}$. Connections between Wasserstein DRO and regularization have been established in several settings (see the details in \citealt{shafieezadeh2019regularization,gao2024wasserstein,yang2022wasserstein}), including for Wasserstein DRCCPs \citep{hota2019data,xie2021distributionally,chen2024data}, and SCA formulation~\eqref{eq_dc} inherits this regularization-based interpretation.
Before presenting the main properties of Algorithm~\ref{alg_DCA}, we refer the interested reader to \Cref{ex:DC_procedure} in Appendix~\ref{appendix_example} for an illustrative example. 

We now establish the convergence of SCA Algorithm~\ref{alg_DCA}, whose iterations are illustrated in Example~\ref{ex:DC_procedure} of Appendix~\ref{appendix_example}. As shown above, every iterate is feasible, i.e., $\bm x^{(k)}\in\Omega(\mathrm{DC})$. Moreover, since $\bm x^{(k)}$ is itself feasible for the $k$-th subproblem~\eqref{eq_dc}, the objective values $\bm c^\top\bm x^{(k)}$ are nonincreasing, and they converge because $\X$ is compact. To show that the iterates also approach critical points of DRCCP~\eqref{eq_dc_orginial}, that is, points $\bar{\bm x}\in\Omega(\mathrm{DC})$ for which there exist $\lambda\geq 0$, $\bm s_1\in\partial g_1(\bar{\bm x})$, and $\bm s_2\in\partial g_2(\bar{\bm x})$ such that $\bm 0\in\bm c+\lambda(\bm s_1-\bm s_2)+\N_\X(\bar{\bm x})$ and $\lambda[g_1(\bar{\bm x})-g_2(\bar{\bm x})]=0$, we impose the following constraint qualification.
\begin{assumption}
\label{assumption_lambda}
For every $\bm x\in\Omega(\mathrm{DC})$ with $g_1(\bm x)=g_2(\bm x)$, we have $\bm 0\notin\partial g_1(\bm x)-\partial g_2(\bm x)+\N_\X(\bm x)$.
\end{assumption}%
Assumption~\ref{assumption_lambda} requires any feasible point \(\bm x\) that is active at the DC constraint not to be a stationary point of the constraint function $g_1-g_2$ over $\X$. Intuitively, the constraint can then be made strictly satisfied by moving within $\X$, which guarantees that every SCA subproblem~\eqref{eq_dc} satisfies Slater's condition (Lemma~\ref{lemma_dc} in Appendix~\ref{appendix_complexity}). When $g_1$ and $g_2$ are differentiable, Assumption~\ref{assumption_lambda} requires that, at every active feasible point, some direction $\bm d$ tangent to $\X$ satisfies $(\nabla g_1(\bm x)-\nabla g_2(\bm x))^\top\bm d<0$. This is Robinson's constraint qualification for $g_1-g_2\le 0$ over $\X$, which is the Mangasarian--Fromovitz constraint qualification when $\X=\Re^n$ \citep{rockafellar1998variational}.
\begin{proposition}\label{prop_converge}
Let $\{\bm x^{(k)}\}_{k \in \Ze_+}$ be the sequence of solutions generated by SCA formulation \eqref{eq_dc} in Algorithm~\ref{alg_DCA} starting from $\bm{x}^{(0)} \in \Omega(\mathrm{DC})$. Then, $\bm x^{(k)} \in \Omega(\mathrm{DC})$ for all \(k\in\Ze_+\), and $\{\bm c^\top \bm x^{(k)}\}_{k \in \Ze_+}$ is a convergent nonincreasing sequence. If, in addition, Assumption~\ref{assumption_lambda} holds, then every cluster point of $\{\bm x^{(k)}\}_{k \in \Ze_+}$ is a critical point of DRCCP~\eqref{eq_dc_orginial}.
\end{proposition}
To further quantify the convergence rate, we denote by $\Delta_k:=\bm c^\top\bm x^{(k)}-\bm c^\top\bm x^{(k+1)}$ the descent in the objective function, and define the criticality gap of $\bar{\bm x}\in\Omega(\mathrm{DC})$ as $\mathcal G(\bar{\bm x}):=\bm c^\top\bar{\bm x}-\sup_{\bm s_2(\bar{\bm x})\in\partial g_2(\bar{\bm x})}\min_{\bm x\in\Omega(\mathrm{DC})_{\bar{\bm x}}}\bm c^\top\bm x$, which, under Assumption~\ref{assumption_lambda}, is zero exactly when $\bar{\bm x}$ is a critical point of DRCCP~\eqref{eq_dc_orginial} (see Appendix~\ref{appendix_complexity}).
\begin{theorem}
\label{thm:iteration-complexity}
Let $\{\bm x^{(k)}\}_{k\in\Ze_+}$ be generated by Algorithm~\ref{alg_DCA} and $\bar M := \max_{\bm x\in \X}\bm c^\top\bm x - \min_{\bm x\in \X}\bm c^\top\bm x$. Then $\min_{1\le i\le k}\mathcal G(\bm x^{(i)})\leq\min_{1\le i\le k}\Delta_i\le \bar M/k$ for every $k\geq 1$. Hence, an $\bar\epsilon$-approximate critical point, i.e., a feasible point $\bar{\bm x}$ with $\mathcal G(\bar{\bm x})\le\bar\epsilon$, is obtained within $\lceil\bar M/\bar\epsilon\rceil+1$ iterations.
\end{theorem}
The proofs of Proposition~\ref{prop_converge} and Theorem~\ref{thm:iteration-complexity} are given in Appendix~\ref{appendix_complexity}.

\section{Conservatism Reduction: Comparison with CVaR Approximation and ALSO-X\#}
\label{sec_comparisons}

In this section, we compare the solutions of SCA Algorithm~\ref{alg_DCA} with CVaR approximation~\eqref{worse_case_cvar_q_formulation_a} and ALSO-X\#~\citep{jiang2025also}. As we illustrated in Example~\ref{example_sensitivity}, both approximations can be overly conservative: the CVaR approximation yields an objective value of approximately three times the optimal value for all $\alpha$, and ALSO-X\# yields the same value for $\alpha\geq 0.76$. 
A natural question is whether SCA, which is designed to mitigate this conservatism, yields better solutions than these two approximations. The answer, however, is not immediate: CVaR approximation~\eqref{worse_case_cvar_q_formulation_a} replaces the nonnegative right-hand side of the constraint $g_1(\bm x)\leq g_2(\bm x)$ in DRCCP~\eqref{eq_dc_orginial} with zero, whereas SCA replaces it with the linearization of $g_2$ at the current iterate, which can be negative. As a result, the SCA subproblem may then exclude the CVaR solution, producing a worse solution.
In Section~\ref{sec_cvar_compare}, we identify initialization rules under which SCA performs no worse than the CVaR approximation, and we show that SCA can return a feasible solution even when the CVaR approximation is infeasible. ALSO-X\#, by contrast, reduces conservatism through bisection on the objective value; in Section~\ref{sec_alsox}, we demonstrate through examples that neither method dominates the other.

\subsection{Outperforming CVaR Approximation After One Additional Iteration}
\label{sec_cvar_compare}

To understand how SCA improves upon the CVaR approximation, we examine the relationship between their respective feasible regions. Specifically, the SCA subproblem contains the feasible region of CVaR approximation \eqref{worse_case_cvar_q_formulation_a} whenever the linearization of \(g_2\) is nonnegative on \(\mathcal X\). %
The following proposition demonstrates that a suitable choice of the initial point and subgradient guarantees this property.

\begin{proposition}
\label{prop_dc_better}
Let $\hat{\bm x}\in\argmin_{\bm x\in\X} g_2(\bm x)$ satisfy $g_1(\hat{\bm x})\leq g_2(\hat{\bm x})$, i.e., $\hat{\bm x}\in\Omega(\mathrm{DC})$. Then there exists $\bm s_2^{(0)}\in\partial g_2(\hat{\bm x})$ with $-\bm s_2^{(0)}\in\N_\X(\hat{\bm x})$. If SCA Algorithm~\ref{alg_DCA} is initialized at $\bm x^{(0)}=\hat{\bm x}$ and uses this subgradient in its first subproblem, then $v^{\mathrm{DC}}\leq v^{\CVaR}$.
\end{proposition}

\begin{proof}
Since $g_2$ is convex and finite on $\Re^n$ and $\hat{\bm x}$ minimizes it over the convex set $\X$, we have $\bm 0\in\partial g_2(\hat{\bm x})+\N_\X(\hat{\bm x})$ \citep[ Theorem~8.15]{rockafellar1998variational}, which yields the desired $\bm s_2^{(0)}$. By the definition of the normal cone, $(\bm s_2^{(0)})^\top(\bm x-\hat{\bm x})\geq 0$ for all $\bm x\in\X$. Together with $g_2\geq 0$, this implies that the linearization of $g_2$ in the first subproblem is nonnegative on $\X$:
\[
g_2(\hat{\bm x})+(\bm s_2^{(0)})^\top(\bm x-\hat{\bm x})\geq 0\quad\text{for all }\bm x\in\X.
\]
Without loss of generality, we assume that \(v^\text{CVaR} < +\infty\), and let $\bm x^{\CVaR}$ be an optimal solution of CVaR approximation~\eqref{worse_case_cvar_q_formulation_a}, so that $g_1(\bm x^{\CVaR})\leq 0$. Consequently,
\[
g_1(\bm x^{\CVaR})\leq 0\leq g_2(\hat{\bm x})+(\bm s_2^{(0)})^\top(\bm x^{\CVaR}-\hat{\bm x}),
\]
from which $\bm x^{\CVaR}$ is feasible for SCA formulation~\eqref{eq_dc} with $k=0$, and hence $\bm c^\top\bm x^{(1)}\leq v^{\CVaR}$. By Proposition~\ref{prop_converge}, the objective values $\bm c^\top\bm x^{(k)}$ are nonincreasing in $k$. Since the output $\bm x^{\mathrm{DC}}$ is an iterate $\bm x^{(k)}$ with $k\geq 1$, we obtain $v^{\mathrm{DC}}\leq\bm c^\top\bm x^{(1)}\leq v^{\CVaR}$.
\QEDA
\end{proof}

Proposition~\ref{prop_dc_better} implies that, with a suitable initialization, a single SCA iteration suffices to match or improve upon CVaR approximation~\eqref{worse_case_cvar_q_formulation_a}. In comparison, bisection-based algorithms offer no comparable guarantee: ALSO-X~\citep{ahmed2017nonanticipative,jiang2022also} refines the CVaR approximation for regular chance constrained programs, and ALSO-X\#~\citep{jiang2025also} extends it to DRCCPs. Although each bisection step solves a convex program of the same type as the CVaR approximation, the number of steps needed to improve on the CVaR solution depends on the initial search interval and the stopping tolerance. \citet{jiang2022also,jiang2025also} report favorable numerical results on their test instances, but no explicit bound on this number is available, so the computational effort needed to outperform the CVaR approximation is instance dependent.

While Proposition~\ref{prop_dc_better} highlights how SCA improves upon an existing CVaR solution, the advantages of our approach extend to cases where the CVaR
approximation \eqref{worse_case_cvar_q_formulation_a} is infeasible, even though the underlying DRCCP remains feasible. This discrepancy arises because the two formulations treat unsafe scenarios differently: the DRCCP truncates the constraint value $g(\bm x,\hat{\rxi})$ at zero and thus limits only the probability of violation, whereas the CVaR approximation also penalizes the magnitude of violations in the upper $\varepsilon$-tail (see Section~\ref{sec_SCA_algorithm}). Consequently, if $g(\bm x,\hat{\rxi})$ stays below a negative threshold with high probability, the resulting safety margin can offset the term $\theta\varepsilon^{-1}\|\bm a_1(\bm x)\|_{p^*}$, no matter how large the magnitude of constraint violation on the remaining scenarios. Such a decision is feasible for the DRCCP, although its large violations may render it infeasible for the CVaR approximation. Proposition~\ref{prop_infeasible} formalizes this argument and shows that SCA initialized at such a decision returns a feasible solution, and Example~\ref{example_cvar_infeasible} in Appendix~\ref{appendix_example} gives a concrete instance.

\begin{subequations}
\begin{proposition}
\label{prop_infeasible}
Consider DRCCP~\eqref{eq_dc_orginial} and its CVaR approximation~\eqref{worse_case_cvar_q_formulation_a}.  Suppose that CVaR approximation~\eqref{worse_case_cvar_q_formulation_a} is infeasible, i.e.,
$v^{\CVaR}=+\infty$. Suppose further that there exist $\hat{\bm x}\in\X$, $\rho>0$, and $\alpha\in[0,\varepsilon)$ such that
\begin{align}
\hat{\Pr}\left\{
\max_{i\in[I]}
\left\{
\bm a_i(\hat{\bm x})^\top\hat{\rxi}-b_i(\hat{\bm x})
\right\}
>-\rho
\right\}
\le \alpha,
\label{eq_margin_prob}
\end{align}
and
\begin{align}
\theta\|\bm a_1(\hat{\bm x})\|_{p^*}
\le
\rho(\varepsilon-\alpha).
\label{eq_margin_radius}
\end{align}
Assume that SCA Algorithm~\ref{alg_DCA} is initialized at 
$\bm x^{(0)}=\hat{\bm x}$. Then, we have $v^{\mathrm{DC}} < v^{\CVaR}$.
\end{proposition}
\end{subequations}
\begin{proof}
We first show that $\hat{\bm x}$ is feasible to DRCCP~\eqref{eq_dc_orginial}. By the definition of CVaR, choosing $\beta=-\rho$ gives
\begin{align*}
\hat{\Pr}\text{-}\CVaR_{1-\varepsilon} \left[ \left( \max_{i\in[I]} \left\{ \bm a_i(\hat{\bm x})^\top\hat{\rxi}-b_i(\hat{\bm x}) \right\} \right)_- \right] \le -\rho + \frac{1}{\varepsilon} \mathbb E_{\hat{\Pr}} \left[ \left( \left( \max_{i\in[I]} \left\{ \bm a_i(\hat{\bm x})^\top\hat{\rxi}-b_i(\hat{\bm x}) \right\} \right)_- +\rho \right)_+ \right].
\end{align*}
According to condition~\eqref{eq_margin_prob}, we have
\begin{align*}
\hat{\Pr}\text{-}\CVaR_{1-\varepsilon}
\left[ \left( \max_{i\in[I]} \left\{ \bm a_i(\hat{\bm x})^\top\hat{\rxi}-b_i(\hat{\bm x}) \right\} \right)_- \right] \le -\rho+\frac{\rho\alpha}{\varepsilon} = -\frac{\rho(\varepsilon-\alpha)}{\varepsilon}.
\end{align*}
Together with condition~\eqref{eq_margin_radius}, this implies
\begin{align*}
\frac{\theta}{\varepsilon}\|\bm a_1(\hat{\bm x})\|_{p^*} + \hat{\Pr}\text{-}\CVaR_{1-\varepsilon} \left[ \left( \max_{i\in[I]} \left\{ \bm a_i(\hat{\bm x})^\top\hat{\rxi}-b_i(\hat{\bm x}) \right\} \right)_- \right] \le 0,
\end{align*}
which further implies that
\begin{align*}
&\frac{\theta}{\varepsilon}\|\bm a_1(\hat{\bm x})\|_{p^*} + \hat{\Pr}\text{-}\CVaR_{1-\varepsilon} \left[ \max_{i\in[I]} \left\{ \bm a_i(\hat{\bm x})^\top\hat{\rxi}-b_i(\hat{\bm x}) \right\} \right]\\
&\qquad\le \hat{\Pr}\text{-}\CVaR_{1-\varepsilon} \left[ \left( \max_{i\in[I]} \left\{ \bm a_i(\hat{\bm x})^\top\hat{\rxi}-b_i(\hat{\bm x}) \right\} \right)_+ \right].
\end{align*}
By the  argument above, $\hat{\bm{x}}$ is feasible to DRCCP~\eqref{eq_dc_orginial}. Therefore, SCA Algorithm~\ref{alg_DCA} can be initialized with a feasible solution $\bm x^{(0)}=\hat{\bm x}$. Since each current iterate is feasible for its SCA subproblem, the subproblems remain feasible, and the objective values are nonincreasing. Together with the compactness of $\X$, these properties ensure that SCA Algorithm~\ref{alg_DCA} terminates for every \(\tau>0\) with a feasible output satisfying $v^{\mathrm{DC}}\le \bm c^\top\hat{\bm x}<+\infty$.
On the other hand, CVaR approximation~\eqref{worse_case_cvar_q_formulation_a} is infeasible by assumption, so $v^{\CVaR}=+\infty$. Thus, $v^{\mathrm{DC}}<+\infty=v^{\CVaR}$. 
\QEDA
\end{proof}

As far as we know, the sufficient conditions in Proposition~\ref{prop_infeasible} are among the first to address the infeasibility issue of CVaR approximation~\eqref{worse_case_cvar_q_formulation_a}. Proposition~\ref{prop_infeasible} provides a sufficient condition for SCA Algorithm~\ref{alg_DCA} to remain effective when the CVaR approximation is infeasible. Condition~\eqref{eq_margin_prob} requires that, under the reference distribution $\hat{\Pr}$, the quantity $\max_{i\in[I]}\{
\bm a_i(\hat{\bm x})^\top\hat{\rxi}-b_i(\hat{\bm x})\}$ exceeds the negative threshold $-\rho$ with probability at most $\alpha$, and condition~\eqref{eq_margin_radius} requires this margin to dominate the nonnegative term $\theta\varepsilon^{-1}\|\bm{a}_1(\hat{\bm x})\|_{p^*}$.

A natural question is whether a similar improvement to that in Proposition~\ref{prop_dc_better} can be obtained under a simpler initialization scheme. The following corollary shows that, under mild conditions satisfied by many DRCCPs in the literature (see, e.g., \citealt{xie2021distributionally,chen2024data}), initializing at the trivial solution $\bm{x}^{(0)}=\bm{0}$ is sufficient to guarantee improvement over the CVaR approximation.

\begin{corollary}
\label{corollary_better}
Suppose that $\bm{a}_i(\bm{0}) = \bm{0}$, $b_i(\bm{0}) \geq 0$ for all $i \in [I]$,  and $ \bm{0}\in\X$. If SCA Algorithm~\ref{alg_DCA} is initialized at $\bm{x}^{(0)} = \bm{0}$ and we choose $ \bm{s}_2^{(0)} = \bm 0$, then its output is better than CVaR approximation~\eqref{worse_case_cvar_q_formulation_a}, i.e., $ v^{\mathrm{DC}}\leq v^{\CVaR}$.
\end{corollary}
\begin{proof}
Since $ \bm{0}\in\X$, we have $\bm 0 \in\argmin_{\bm x\in\X}\hat{\Pr}\text{-}\CVaR_{1-\varepsilon}[ \max_{i\in[I]} \{\bm{a}_i(\bm{x})^\top\hat{\rxi}-b_i(\bm{x})\}_+]$. With the assumptions $\bm{a}_i(\bm{0}) = \bm{0}$, $b_i(\bm{0}) \geq 0$ for all $i \in [I]$, we always have 
\begin{align*}
\theta\varepsilon^{-1}\left\|\bm{a}_1(\bm 0)\right\|_{p^*} + \hat{\Pr}\text{-}\CVaR_{1-\varepsilon}\left[ \max_{i\in[I]} \left\{\bm{a}_i(\bm 0)^\top\hat{\rxi}-b_i(\bm 0)\right\}\right] \leq \hat{\Pr}\text{-}\CVaR_{1-\varepsilon}\left[ \max_{i\in[I]} \left\{\bm{a}_i(\bm 0)^\top\hat{\rxi}-b_i(\bm 0)\right\}_+\right].
\end{align*}
Following the results in Proposition~\ref{prop_dc_better}, we obtain $ v^{\mathrm{DC}}\leq v^{\CVaR}$.
\QEDA
\end{proof}

The premises of Corollary~\ref{corollary_better} are satisfied by packing problems, such as the distributionally robust multidimensional knapsack problems~\citep{xie2021distributionally, ji2021data}, where $\X=[0,1]^n$ and each uncertain constraint takes the form $\tilde{\bm\xi}_i^\top\bm x\leq B_i$ with capacity $B_i>0$, so that $\bm a_i(\bm 0)=\bm 0$ and $b_i(\bm 0)=B_i>0$.

Even when the conditions in Proposition~\ref{prop_dc_better} or Corollary~\ref{corollary_better} are satisfied, it may still occur that SCA Algorithm~\ref{alg_DCA} terminates at the same objective value as CVaR approximation~\eqref{worse_case_cvar_q_formulation_a}, i.e., $v^{\mathrm{DC}}=v^{\CVaR}$. This situation can arise when $\bm{x}^\CVaR$, an optimal solution of CVaR approximation~\eqref{worse_case_cvar_q_formulation_a}, also minimizes $\hat{\Pr}\text{-}\CVaR_{1-\varepsilon}[ \max_{i\in[I]} \{\bm{a}_i(\bm{x})^\top\hat{\rxi}-b_i(\bm{x})\}_+]$. We refer the interested reader to \Cref{example_revisit_SCA} in Appendix~\ref{appendix_example} for an illustration.

\subsection{SCA Algorithm~\ref{alg_DCA} and ALSO-X\# Are Not Directly Comparable}
\label{sec_alsox}
ALSO-X\#~\citep{jiang2025also} refines the CVaR approximation through a bisection search on the objective value, as summarized in Algorithm~\ref{alg_alsox_sharp}. For each candidate level $t$, the lower level minimizes the CVaR constraint function $g_1$ over $\{\bm x\in\X:\bm c^\top\bm x\leq t\}$, as in~\eqref{eq_also_sharp_ori}, and the upper level accepts $t$ as an upper bound on $v^*$ if the resulting minimizer $\bm x^{\as}$ satisfies the exact DRCCP constraint~\eqref{alsoxs_drccp_formualtion_upper}.
\begin{subequations}
\begin{algorithm}[H]
\caption{{ALSO-X\#} \citep{jiang2025also}}
\label{alg_alsox_sharp}
\begin{algorithmic}[1]
\State \textbf{Input:} Let $\delta_1$ denote the stopping tolerance parameter, $t_L$ and $t_U$ be the known lower and upper bounds of the optimal value of DRCCP~\eqref{1_wasser_reformualtion}, respectively 
\While {$t_U-t_L>\delta_1$}
\State {Let $t=(t_L+t_U)/2$ and $\bm x^{\as}$ be an optimal solution of the lower-level ALSO-X\# \eqref{eq_also_sharp_ori}:
\begin{align}
\bm x^{\as}\in\argmin _{\bm {x}\in \X}\left\{ \theta\varepsilon^{-1}\left\|\bm{a}_1(\bm{x})\right\|_{p^*} + \hat{\Pr}\text{-}\CVaR_{1-\varepsilon}\left[ \max_{i\in[I]} \left\{\bm{a}_i(\bm{x})^\top\hat{\rxi}-b_i(\bm{x})\right\}\right]\colon \bm{c}^\top \bm{x} \leq t \right\}\label{eq_also_sharp_ori}
\end{align}
\State Let $t_U=t$ if $\bm x^{\as}$ satisfies the checking condition of the upper-level ALSO-X\# \eqref{alsoxs_drccp_formualtion_upper}:
\begin{align}
\theta\varepsilon^{-1}\left\|\bm{a}_1(\bm x^{\as})\right\|_{p^*}+\hat{\Pr}\text{-}\CVaR_{1-\varepsilon}\left[-\min_{i\in [I]}\left({b}_i(\bm x^{\as})-\bm a_i(\bm x^{\as})^\top\hat{\rxi}\right)_+\right]\leq 0; \label{alsoxs_drccp_formualtion_upper}
\end{align}
\,\quad otherwise, $t_L=t$}
\EndWhile
\State \textbf{Output:} A solution $\bm x^{\as}$ and its objective value $\bm c^\top \bm x^{\as}$ to DRCCP~\eqref{1_wasser_reformualtion} \end{algorithmic}
\end{algorithm}
\end{subequations}
ALSO-X\# and SCA use the exact DRCCP constraint at different stages, thereby affecting how each method behaves numerically: ALSO-X\# generates its candidates by minimizing $g_1$ and uses the exact DRCCP constraint only to check their feasibility. By~\eqref{CVaR_comonotonicity}, $g_1$ equals the left-hand side of DRCCP~\eqref{1_wasser_reformualtion} plus the term $\hat{\Pr}\text{-}\CVaR_{1-\varepsilon}[g(\bm x,\hat{\rxi})_+]$, which penalizes the magnitude of constraint violations that the DRCCP ignores (Region~III in Figure~\ref{fig:tail_decomposition}). Minimizing $g_1$ therefore steers the candidates away from points that severely violate the constraint under a few scenarios, even though such points can be feasible for the DRCCP. When this occurs at a level $t$, the minimizer of $g_1$ is not feasible for the DRCCP, and ALSO-X\# rejects $t$ even though another point at that level satisfies the DRCCP constraint.
SCA, in contrast, retains the DC structure of the exact constraint in every subproblem and offsets the Region~III term through the linearization of $g_2$. However, each subproblem is restricted to an inner approximation of the feasible region around the current iterate, and the algorithm stops once this approximation yields no further improvement.

The examples in Appendix~\ref{appendix_example} show that either effect can dominate. In Example~\ref{example_sensitivity}, the conservatism stems from scenarios in Region~III, including scenario $\hat{\rxi}^2$ perturbed through the parameter $\alpha\in[0.7, 0.9]$; for $\alpha\geq 0.76$, ALSO-X\# returns the CVaR value, whereas SCA attains the value $10/9$ in two iterations (\Cref{example_sensitivity_detail}). SCA also outperforms ALSO-X\# in \Cref{example_compare_also_a}, with $-0.742$ versus $-0.522$. In \Cref{example_compare_also_b}, by contrast, SCA stops at $-0.939$, whereas ALSO-X\#, whose search is not tied to a current iterate, reaches $-1.049$. Hence, neither method dominates the other. Their computational costs also differ: ALSO-X\# always performs $\lceil\log_2((t_U-t_L)/\delta_1)\rceil$ bisection steps, which equals $13$ in both examples, whereas SCA converges in three iterations.

The two methods can also be combined. Since the output of ALSO-X\# satisfies the DRCCP constraint, it can serve as the initial point of SCA Algorithm~\ref{alg_DCA}, and the monotonicity of SCA then guarantees $v^{\mathrm{DC}}\leq\bm c^\top\bm x^{\as}$. The combined procedure therefore performs no worse than ALSO-X\#.

\section{Numerical Results}
\label{sec_Numerical_study}
In this section, we present numerical evidence on the effectiveness of the proposed methods. All instances are implemented in Python 3.11, with calls to the Gurobi solver (version 12.0.2, default settings), on a personal computer with a 12th-Gen Intel(R) Core(TM) i5-12500H processor and 16GB of memory. 
To assess the performance of the different algorithms, we report relative improvement in objective value. For any pairwise comparison, we define
\begin{align*}
\textrm{Improvement }(\%)=
\frac{\textrm{reference value}-\textrm{algorithm value}}
{|\textrm{reference value}|}\times 100\%.
\end{align*}
We use $``\textrm{UB}" $ and $``\textrm{LB}" $ to denote the best upper bound and the best lower bound found by the solver. Since some instances may not be solved to optimality within the time limit, we use ``GAP" to denote the optimality gap as $\textrm{GAP } (\%)=(| \textrm{UB} -\textrm{LB} |)/(|\textrm{LB}|)\times 100$.  In all our experiments, we set the time limit for each instance to 3600 seconds and the default optimality gap tolerance to 0.01\%.  

\subsection{Portfolio Optimization: An Illustrative Example Adapted from \texorpdfstring{\citet{ahmed2018relaxations}}{Ahmed and Xie (2018)}}
\label{sec_num_single}
In this subsection, we consider the portfolio optimization problem studied in \citet{ahmed2018relaxations}, with related formulations also appearing in \citet{qiu2014covering,chen2024data,ho2023strong}. In this problem, the decision-maker seeks a minimum-cost investment vector $\x\in\X$ such that the portfolio return exceeds a prescribed target $w$ with high confidence. To exclude the degenerate solution $\bm 0$, we assume that $\X=\{\bm x\in\Re^n_+:\sum_{i\in[n]}x_i\geq 0.01, 0\leq x_i\leq 10000, \forall i\in[n]\}$. Let $\tilde{\rxi}\in\Re^n$ denote the random return vector of the $n$ assets. The problem can be written as an instance of DRCCP~\eqref{eq_drccp_org_ch2}:
\begin{align}
v^*
=\min_{\x\in\X}\left\{
\bm c^\top \x:\ 
\inf_{\Pr\in\P}
\Pr\left\{\trxi^{\top}\x \ge w\right\} \ge  1-\varepsilon
\right\}.
\label{eq:portfolio_drccp}
\end{align}
We generate instances with $n=50$ assets and target return level $w=1$, and we sample the cost coefficients $c_1,\dots,c_{50}$ independently and uniformly from $\{1,\dots,100\}$. The true (nominal) return vector $\trxi$ has i.i.d.\ components that are uniform on $[0.8,1.5]$. To reflect abnormal market conditions, the observed historical returns are drawn not from $\trxi$ but from a contaminated distribution: with probability $0.92$, an observation is a draw of $\trxi$, and with probability $0.08$, it is a deteriorated draw $\trxi_o$, obtained by multiplying each component of a draw of $\trxi$ by an independent factor uniform on $[0.02,0.20]$. We draw $N=100$ observations from this contaminated distribution, denoted by $\hat{\rxi}^1,\dots,\hat{\rxi}^{100}$, and let $\hat{\Pr}$ denote the associated empirical distribution. In our implementation, the contamination is realized by drawing $100$ samples of $\trxi$ and scaling a randomly selected $8\%$ of them in this way.

We compare four approaches: CVaR approximation~\eqref{worse_case_cvar_q_formulation_a}, ALSO-X\# Algorithm~\ref{alg_alsox_sharp}, our SCA Algorithm~\ref{alg_DCA}, and the mixed-integer benchmark solved via a Big-M reformulation (see, e.g., \citealt{xie2021distributionally,chen2024data}). 
We use the $2$-norm Wasserstein ambiguity set and consider the Wasserstein radii $\theta\in\{0.05,0.10,0.20\}$ and  the risk levels $\varepsilon\in\{0.10,0.20\}$ in~\eqref{eq:portfolio_drccp}. In SCA Algorithm~\ref{alg_DCA}, we choose $\tau$ as $10^{-3}$ and set the initial solution as $ x^{(0)}_i=10000$ for each $i\in[n]$. At each feasible iterate $\bm x^{(k)}$, we select $\bm s_2^{(k)}=-(N\varepsilon)^{-1}
\sum_{j=1}^{N}\mathbf 1_{\{(\hat{\rxi}^{j})^\top\bm x^{(k)}<w\}}\hat{\rxi}^{j}\in\partial g_2(\bm x^{(k)})$. In ALSO-X\# Algorithm~\ref{alg_alsox_sharp}, we choose $t_L$ as $0$, $t_U$ as $2000$, and $\delta_1 = 10^{-3}$. We report the average performance over five random instances. In this setting, the Big-M method solves all instances to optimality very quickly, confirming that the benchmark solutions used for comparison are exact. Table~\ref{tab:seed_avg} reports the average objective values $v^{\CVaR}$, $v^{\as}$, $v^{\mathrm{DC}}$, and $v^*$, where $v^*=v^{\mathrm{BigM}}$ is obtained by the Big-M method. For each method, let $\mathrm{rt}^{\CVaR}$, $\mathrm{rt}^{\as}$, $\mathrm{rt}^{\mathrm{DC}}$, and $\mathrm{rt}^{\mathrm{BigM}}$ denote the average running time (in seconds) across these five random instances. The column ``Imp$_1$\ (\%)'' denotes the average relative objective reduction achieved by SCA Algorithm~\ref{alg_DCA} with respect to CVaR approximation~\eqref{worse_case_cvar_q_formulation_a} over five instances, whereas ``Imp$_2$ (\%)'' denotes the average relative objective reduction of SCA Algorithm~\ref{alg_DCA} with respect to ALSO-X\# Algorithm~\ref{alg_alsox_sharp}. Overall, SCA Algorithm~\ref{alg_DCA} consistently improves upon both CVaR approximation~\eqref{worse_case_cvar_q_formulation_a} and ALSO-X\# Algorithm~\ref{alg_alsox_sharp}. Compared with the CVaR approximation, the improvement is particularly pronounced for smaller $\theta$, where Imp$_1$ ranges from about $30\%$ to $50\%$. Compared with ALSO-X\#, it also attains lower objective values, although the gap becomes small in some more difficult instances, especially when $\varepsilon = 0.20$. Moreover, the objective values obtained by the SCA algorithm are nearly identical to those of the exact mixed-integer benchmark 
across all tested instances. From a computational perspective, SCA Algorithm~\ref{alg_DCA} requires only $0.03$--$0.05$ seconds on average, which is comparable to CVaR approximation~\eqref{worse_case_cvar_q_formulation_a} and faster than both ALSO-X\# Algorithm~\ref{alg_alsox_sharp} and the Big-M method.

\begin{table}[t]
\centering
\caption{Numerical results of portfolio optimization. Here $v^{\CVaR}$, $v^{\as}$, $v^{\mathrm{DC}}$, and $v^*$ denote the average objective values over five random instances of CVaR approximation ~\eqref{worse_case_cvar_q_formulation_a}, ALSO-X\# Algorithm~\ref{alg_alsox_sharp}, SCA Algorithm~\ref{alg_DCA}, and the exact Big-M method, respectively.}
\label{tab:seed_avg}
\small
\setlength{\tabcolsep}{4pt}
\begin{tabular}{ccrrrrrrrrrr}
\toprule
$\varepsilon$ & $\theta$
& $v^{\CVaR}$
& $v^{\as}$
& $v^{\mathrm{DC}}$
& $v^*$
& $\mathrm{rt}^{\CVaR}$
& $\mathrm{rt}^{\as}$
& $\mathrm{rt}^{\mathrm{DC}}$
& $\mathrm{rt}^{\mathrm{BigM}}$
& Imp$_1$ (\%)
& Imp$_2$ (\%) \\
\midrule
0.10 & 0.05 &  89.836 &  60.040 &  43.711 &  43.711 & 0.01 & 0.21 & 0.04 & 0.55 & 51.03 & 28.27 \\
0.10 & 0.10 & 218.700 & 165.720 & 125.067 & 125.067 & 0.01 & 0.19 & 0.05 & 1.10 & 42.65 & 25.70 \\
0.10 & 0.20 &1207.776 &1134.540 &1133.908 &1133.908 & 0.01 & 0.23 & 0.03 & 1.35 &  6.32 &  0.06 \\
\midrule
0.20 & 0.05 &  11.449 &   9.810 &   7.793 &   7.793 & 0.01 & 0.26 & 0.03 & 1.80 & 31.34 & 20.56 \\
0.20 & 0.10 &  17.890 &  12.640 &  12.544 &  12.544 & 0.01 & 0.24 & 0.03 & 4.50 & 29.12 &  0.76 \\
0.20 & 0.20 &  36.478 &  26.340 &  26.050 &  26.050 & 0.01 & 0.22 & 0.03 & 7.93 & 28.41 &  1.10 \\
\midrule
\multicolumn{6}{r}{Average}
& 0.01 & 0.23 & 0.03 & 2.87 & 31.48 & 12.74 \\
\bottomrule
\end{tabular}
\end{table}

\subsection{Joint Chance Constraints in Portfolio Optimization}
\label{sec_num_joint}
We next consider a joint chance constraint version of the portfolio optimization problem in Section~\ref{sec_num_single}. In contrast to the individual chance constraint~\eqref{eq:portfolio_drccp} considered in Section~\ref{sec_num_single}, we now require multiple return targets to be satisfied simultaneously with high probability.
Given $K$ target return levels $w_1,\dots,w_K$, the decision-maker seeks a minimum-cost portfolio $\bm x\in\X$ such that all return requirements are satisfied jointly with probability at least $1-\varepsilon$, where $\X=\{\bm x\in\Re^n_+:\sum_{i\in[n]}x_i\geq 0.01, 0\leq x_i\leq 10000, \forall i\in[n]\}$. This problem can be written as an instance of DRCCP~\eqref{eq_drccp_org_ch2}:
\begin{align}
v^*= \min_{\x\in\X} \left\{\bm c^\top \x:\inf_{\Pr\in\P}\Pr\left\{\min_{k\in[K]}\left\{\tilde{\rxi}_k^\top\x - w_k \right\}\geq 0\right\}\ge 1-\varepsilon\right\}.
\label{eq:portfolio_joint_drccp}
\end{align}
For the numerical experiments, we generate instances with $n=50$ assets and $K=5$ joint return constraints, following the setup of Section~\ref{sec_num_single} with three changes: the target levels $w_1,\dots,w_K$ are sampled independently and uniformly from $[0.8,1.2]$; the true return vector $\trxi=(\trxi_1,\dots,\trxi_K)\in\Re^{nK}$ has i.i.d.\ components that are uniform on $[0.8,1.5]$, and contamination scales all of its components; and we draw $N=1000$ observations, denoted by $\hat{\rxi}^1,\dots,\hat{\rxi}^{1000}$, with associated empirical distribution $\hat{\Pr}$.

We compare the same four approaches, using the same $2$-norm Wasserstein ambiguity set and parameter settings as in Section~\ref{sec_num_single}. Because of the larger number of scenarios ($N=1000$), the Big-M formulation is computationally demanding: none of the instances is solved to proven optimality within the time limit, so $v^{\mathrm{BigM}}$ in Table~\ref{tab:joint_seed_avg} denotes the best upper bound found. Table~\ref{tab:joint_seed_avg} reports the average results over five random instances, with the same columns as Table~\ref{tab:seed_avg} plus ``Gap (\%)'', the average optimality gap of the Big-M method.

\begin{table}[t]
\centering
\caption{Numerical results of portfolio optimization for the joint chance constraint case. Here $v^{\CVaR}$, $v^{\as}$, $v^{\mathrm{DC}}$, and $v^{\mathrm{BigM}}$ denote the average objective values over five random instances of CVaR approximation~\eqref{worse_case_cvar_q_formulation_a}, ALSO-X\# Algorithm~\ref{alg_alsox_sharp}, SCA Algorithm~\ref{alg_DCA}, and the best upper bound obtained by the Big-M method, respectively.}
\label{tab:joint_seed_avg}
\small
\setlength{\tabcolsep}{2.15pt}
\begin{tabular}{ccrrrrrrrrrrr}
\toprule
$\varepsilon$ & $\theta$
& $v^{\CVaR}$
& $v^{\as}$
& $v^{\mathrm{DC}}$
& $v^{\mathrm{BigM}}$
& $\mathrm{rt}^{\CVaR}$
& $\mathrm{rt}^{\as}$
& $\mathrm{rt}^{\mathrm{DC}}$
& $\mathrm{rt}^{\mathrm{BigM}}$
& Gap (\%)
& Imp$_1$ (\%)
& Imp$_2$ (\%) \\
\midrule
0.10 & 0.05 &   65.283 &   54.125 &   39.442 &   39.584 & 0.73 &  7.70 & 3.33 & 3600.00 & 17.58 & 39.58 & 27.13 \\
0.10 & 0.10 &  235.172 &  192.903 &  143.291 &  143.291 & 1.31 &  7.36 & 6.05 & 3600.00 & 17.57 & 39.07 & 25.72 \\
0.10 & 0.20 & 1611.734 & 1584.727 & 1584.586 & 1604.490 & 1.14 &  6.68 & 4.20 & 3600.00 & 59.10 &  1.68 &  0.01 \\
\midrule
0.20 & 0.05 &    5.786 &    3.965 &    3.958 &    3.958 & 1.44 & 11.46 & 4.65 & 3600.00 & 21.42 & 31.60 &  0.18 \\
0.20 & 0.10 &    9.520 &    6.679 &    6.652 &    6.666 & 1.47 &  8.80 & 4.72 & 3600.00 & 26.88 & 30.12 &  0.40 \\
0.20 & 0.20 &   24.149 &   18.207 &   17.828 &   17.828 & 1.61 &  7.19 & 4.34 & 3600.00 & 27.23 & 26.18 &  2.08 \\
\midrule
\multicolumn{6}{r}{Average}
& 1.28 & 8.20 & 4.55 & 3600.00 & 28.30 & 28.04 & 9.25 \\
\bottomrule
\end{tabular}
\end{table}

As shown in Table~\ref{tab:joint_seed_avg}, SCA Algorithm~\ref{alg_DCA} consistently outperforms both CVaR approximation~\eqref{worse_case_cvar_q_formulation_a} and ALSO-X\# Algorithm~\ref{alg_alsox_sharp} in terms of objective value. The relative improvement over the CVaR approximation ranges from about $26\%$ to $40\%$ in most instances, and the improvement over ALSO-X\# is positive in all instances. 
These results show that the SCA algorithm has a clear advantage in solution quality.
In terms of computation time, SCA Algorithm~\ref{alg_DCA} is much more efficient than the Big-M formulation, which reaches the time limit in all instances. As a result, the Big-M method returns incumbent solutions with noticeable optimality gaps. Moreover, the objective values produced by SCA Algorithm~\ref{alg_DCA} are generally very close to those obtained by the Big-M formulation, and in some instances even better than the incumbent values returned when solving the Big-M formulation within the time limit.

One case with $(\varepsilon,\theta)=(0.10,0.20)$ deserves additional discussion, where the relative improvement of SCA Algorithm~\ref{alg_DCA} over CVaR approximation~\eqref{worse_case_cvar_q_formulation_a} is only about $1.68\%$, and the improvement over ALSO-X\# Algorithm~\ref{alg_alsox_sharp} is nearly zero. These results should not be interpreted as evidence that the CVaR approximation is less conservative or that the three methods perform similarly in this case. Instead, this instance appears to be particularly difficult, as shown by the fact that the Big-M formulation still has an average optimality gap of about $59\%$ after reaching the time limit. Since even the benchmark formulation remains far from convergence, the small differences among the three approximation methods are more likely due to the difficulty of the instance itself than to any reduction in the conservatism of CVaR-based models. Even so, the SCA algorithm retains its computational advantage.

Overall, these results show that SCA Algorithm~\ref{alg_DCA} provides a good balance between solution quality and computational efficiency. It consistently gives strong objective values, improves upon both CVaR approximation~\eqref{worse_case_cvar_q_formulation_a} and ALSO-X\# Algorithm~\ref{alg_alsox_sharp}, and does so at a much lower computational cost than the Big-M method.

\newpage
\bibliography{project}

@article{hong2011sequential,
	title={Sequential convex approximations to joint chance constrained programs: A Monte Carlo approach},
	author={Hong, L Jeff and Yang, Yi and Zhang, Liwei},
	journal={Operations Research},
	volume={59},
	number={3},
	pages={617--630},
	year={2011},
	publisher={INFORMS}
}

@techreport{slater1950lagrange,
  author = {Slater, Morton},
  title = {Lagrange Multipliers Revisited},
  institution = {Cowles Commission for Research in Economics},
  type = {Cowles Commission Discussion Paper: Mathematics},
  number = {403},
  year = {1950}
}

@article{dhaene2002concept,
  title={The concept of comonotonicity in actuarial science and finance: theory},
  author={Dhaene, Jan and Denuit, Michel and Goovaerts, Marc J and Kaas, Rob and Vyncke, David},
  journal={Insurance: Mathematics and Economics},
  volume={31},
  number={1},
  pages={3--33},
  year={2002},
  publisher={Elsevier}
}

@article{zhang2024joint,
  title={Joint chance-constrained program based electric vehicles optimal dispatching strategy considering drivers' response uncertainty},
  author={Zhang, Kaizhe and Xu, Yinliang and Sun, Hongbin},
  journal={Applied Energy},
  volume={356},
  pages={122313},
  year={2024},
  publisher={Elsevier}
}

@article{hanasusanto2017ambiguous,
  title={Ambiguous joint chance constraints under mean and dispersion information},
  author={Hanasusanto, Grani A and Roitch, Vladimir and Kuhn, Daniel and Wiesemann, Wolfram},
  journal={Operations Research},
  volume={65},
  number={3},
  pages={751--767},
  year={2017},
  publisher={Informs}
}

@article{jiang2025terminator,
  title={The terminator: An integration of inner and outer approximations for solving Wasserstein distributionally robust chance constrained programs via variable fixing},
  author={Jiang, Nan and Xie, Weijun},
  journal={INFORMS Journal on Computing},
  volume={37},
  number={2},
  pages={381--412},
  year={2025},
  publisher={INFORMS}
}

@book{denneberg1994non,
  title={Non-additive measure and integral},
  author={Denneberg, Dieter},
  volume={27},
  year={1994},
  publisher={Springer Science \& Business Media}
}

@article{shen2025convex,
  title={Convex chance-constrained programs with Wasserstein ambiguity},
  author={Shen, Haoming and Jiang, Ruiwei},
  journal={Operations Research},
  volume={73},
  number={4},
  pages={2264--2280},
  year={2025},
  publisher={INFORMS}
}

@article{ho2023strong,
  title={Strong formulations for distributionally robust chance-constrained programs with left-hand side uncertainty under Wasserstein ambiguity},
  author={Ho-Nguyen, Nam and Kilin{\c{c}}-Karzan, Fatma and K{\"u}{\c{c}}{\"u}kyavuz, Simge and Lee, Dabeen},
  journal={INFORMS Journal on Optimization},
  volume={5},
  number={2},
  pages={211--232},
  year={2023},
  publisher={INFORMS}
}

@article{wang2026proximal,
  title={A proximal difference-of-convex algorithm for sample average approximation of chance constrained programming},
  author={Wang, Peng and Jiang, Rujun and Kong, Qingyuan and Balzano, Laura},
  journal={INFORMS Journal on Computing},
  volume={38},
  number={1},
  pages={315--339},
  year={2026},
  publisher={INFORMS}
}

@article{ji2021data,
  title={Data-driven distributionally robust chance-constrained optimization with Wasserstein metric},
  author={Ji, Ran and Lejeune, Miguel A},
  journal={Journal of Global Optimization},
  volume={79},
  number={4},
  pages={779--811},
  year={2021},
  publisher={Springer}
}

@incollection{yang2022wasserstein,
  title={Wasserstein regularization for 0--1 loss},
  author={Yang, Zhen and Gao, Rui},
  journal={Optimization Online Preprint},
  year={2022}
}

@article{gao2024wasserstein,
  title={Wasserstein distributionally robust optimization and variation regularization},
  author={Gao, Rui and Chen, Xi and Kleywegt, Anton J},
  journal={Operations Research},
  volume={72},
  number={3},
  pages={1177--1191},
  year={2024},
  publisher={Informs}
}

@article{rockafellar2000optimization,
  title={Optimization of conditional value-at-risk},
  author={Rockafellar, R Tyrrell and Uryasev, Stanislav},
  journal={Journal of risk},
  volume={2},
  number={3},
  pages={21--41},
  year={2000}
}

@article{zymler2013distributionally,
  title={Distributionally robust joint chance constraints with second-order moment information},
  author={Zymler, Steve and Kuhn, Daniel and Rustem, Ber{\c{c}}},
  journal={Mathematical Programming},
  volume={137},
  pages={167--198},
  year={2013},
  publisher={Springer}
}

@article{kuccukyavuz2022chance,
  title={Chance-constrained optimization under limited distributional information: A review of reformulations based on sampling and distributional robustness},
  author={K{\"u}{\c{c}}{\"u}kyavuz, Simge and Jiang, Ruiwei},
  journal={EURO Journal on Computational Optimization},
  volume={10},
  pages={100030},
  year={2022},
  publisher={Elsevier}
}

@article{chen2023approximations,
  title={On approximations of data-driven chance constrained programs over Wasserstein balls},
  author={Chen, Zhi and Kuhn, Daniel and Wiesemann, Wolfram},
  journal={Operations Research Letters},
  volume={51},
  number={3},
  pages={226--233},
  year={2023},
  publisher={Elsevier}
}

@article{shen2023chance,
	title={Chance-constrained set covering with Wasserstein ambiguity},
	author={Shen, Haoming and Jiang, Ruiwei},
	journal={Mathematical Programming},
	volume={198},
	number={1},
	pages={621--674},
	year={2023},
	publisher={Springer}
}

@article{xie2020bicriteria,
	title={Bicriteria approximation of chance-constrained covering problems},
	author={Xie, Weijun and Ahmed, Shabbir},
	journal={Operations Research},
	volume={68},
	number={2},
	pages={516--533},
	year={2020},
	publisher={INFORMS}
}

@article{qiu2014covering,
	title={Covering linear programming with violations},
	author={Qiu, Feng and Ahmed, Shabbir and Dey, Santanu S and Wolsey, Laurence A},
	journal={INFORMS Journal on Computing},
	volume={26},
	number={3},
	pages={531--546},
	year={2014},
	publisher={INFORMS}
}

@article{shafieezadeh2019regularization,
  title={Regularization via mass transportation},
  author={Shafieezadeh-Abadeh, Soroosh and Kuhn, Daniel and Esfahani, Peyman Mohajerin},
  journal={Journal of Machine Learning Research},
  volume={20},
  number={103},
  pages={1--68},
  year={2019}
}

@article{hong2014conditional,
  title={Conditional value-at-risk approximation to value-at-risk constrained programs: A remedy via Monte Carlo},
  author={Hong, L Jeff and Hu, Zhaolin and Zhang, Liwei},
  journal={INFORMS Journal on Computing},
  volume={26},
  number={2},
  pages={385--400},
  year={2014},
  publisher={INFORMS}
}

@article{shi2014optimal,
  title={Optimal stochastic coordinated beamforming for wireless cooperative networks with CSI uncertainty},
  author={Shi, Yuanming and Zhang, Jun and Letaief, Khaled B},
  journal={IEEE Transactions on Signal Processing},
  volume={63},
  number={4},
  pages={960--973},
  year={2015},
  publisher={IEEE}
}

@article{tao1997convex,
  title={Convex analysis approach to DC programming: theory, algorithms and applications},
  author={Tao, Pham Dinh and An, LT Hoai},
  journal={Acta mathematica vietnamica},
  volume={22},
  number={1},
  pages={289--355},
  year={1997}
}

@article{blanchet2019quantifying,
	title={Quantifying distributional model risk via optimal transport},
	author={Blanchet, Jose and Murthy, Karthyek},
	journal={Mathematics of Operations Research},
	volume={44},
	number={2},
	pages={565--600},
	year={2019},
	publisher={INFORMS}
}

@article{ho2022distributionally,
  title={Distributionally robust chance-constrained programs with right-hand side uncertainty under Wasserstein ambiguity},
  author={Ho-Nguyen, Nam and K{\i}l{\i}n{\c{c}}-Karzan, Fatma and K{\"u}{\c{c}}{\"u}kyavuz, Simge and Lee, Dabeen},
  journal={Mathematical Programming},
  volume={196},
  number={1},
  pages={641--672},
  year={2022},
  publisher={Springer}
}

@article{gao2023distributionally,
	title={Distributionally robust stochastic optimization with Wasserstein distance},
	author={Gao, Rui and Kleywegt, Anton},
	journal={Mathematics of Operations Research},
	volume={48},
	number={2},
	pages={603--655},
	year={2023},
	publisher={INFORMS}
}

@article{nemirovski2007convex,
  title={Convex approximations of chance constrained programs},
  author={Nemirovski, Arkadi and Shapiro, Alexander},
  journal={SIAM Journal on Optimization},
  volume={17},
  number={4},
  pages={969--996},
  year={2007},
  publisher={SIAM}
}

@article{chen2024data,
  title={Data-driven chance constrained programs over Wasserstein balls},
  author={Chen, Zhi and Kuhn, Daniel and Wiesemann, Wolfram},
  journal={Operations Research},
  volume={72},
  number={1},
  pages={410--424},
  year={2024},
  publisher={INFORMS}
}

@article{xie2021distributionally,
  title={On distributionally robust chance constrained programs with Wasserstein distance},
  author={Xie, Weijun},
  journal={Mathematical Programming},
  volume={186},
  number={1},
  pages={115--155},
  year={2021},
  publisher={Springer}
}

@article{xie2017distributionally,
  title={Distributionally robust chance constrained optimal power flow with renewables: A conic reformulation},
  author={Xie, Weijun and Ahmed, Shabbir},
  journal={IEEE Transactions on Power Systems},
  volume={33},
  number={2},
  pages={1860--1867},
  year={2018},
  publisher={IEEE}
}

@article{luedtke2010integer,
  title={An integer programming approach for linear programs with probabilistic constraints},
  author={Luedtke, James and Ahmed, Shabbir and Nemhauser, George L},
  journal={Mathematical programming},
  volume={122},
  number={2},
  pages={247--272},
  year={2010},
  publisher={Springer}
}

@article{ahmed2018relaxations,
  title={Relaxations and approximations of chance constraints under finite distributions},
  author={Ahmed, Shabbir and Xie, Weijun},
  journal={Mathematical Programming},
  volume={170},
  pages={43--65},
  year={2018},
  publisher={Springer}
}

@article{ahmed2017nonanticipative,
  title={Nonanticipative duality, relaxations, and formulations for chance-constrained stochastic programs},
  author={Ahmed, Shabbir and Luedtke, James and Song, Yongjia and Xie, Weijun},
  journal={Mathematical Programming},
  volume={162},
  number={1},
  pages={51--81},
  year={2017},
  publisher={Springer}
}

@article{jiang2022also,
	title={ALSO-X and ALSO-X+: Better Convex Approximations for Chance Constrained Programs},
	author={Jiang, Nan and Xie, Weijun},
	journal={Operations Research},
	volume = {70},
	number = {6},
	pages = {3581-3600},
	year={2022},
	publisher={INFORMS}
}

@incollection{li2026solving,
  title={Solving Chance Constrained Programs via a Penalty based Difference of Convex Approach},
  author={Li, Zhiping and Jiang, Nan and Jiang, Rujun},
  journal={Optimization Online Preprint},
  year={2026}
}

@article{jiang2025also,
  title={ALSO-X\#: Better convex approximations for distributionally robust chance constrained programs},
  author={Jiang, Nan and Xie, Weijun},
  journal={Mathematical Programming},
  volume={213},
  pages={575–638},
  year={2025},
  publisher={Springer}
}

@article{liu2023risk,
  title={Risk-based robust statistical learning by stochastic difference-of-convex value-function optimization},
  author={Liu, Junyi and Pang, Jong-Shi},
  journal={Operations Research},
  volume={71},
  number={2},
  pages={397--414},
  year={2023},
  publisher={INFORMS}
}

@article{mohajerin2018data,
  title={Data-driven distributionally robust optimization using the Wasserstein metric: Performance guarantees and tractable reformulations},
  author={Mohajerin Esfahani, Peyman and Kuhn, Daniel},
  journal={Mathematical Programming},
  volume={171},
  number={1},
  pages={115--166},
  year={2018},
  publisher={Springer}
}

@inproceedings{hota2019data,
  title={Data-driven chance constrained optimization under Wasserstein ambiguity sets},
  author={Hota, Ashish R and Cherukuri, Ashish and Lygeros, John},
  booktitle={2019 American Control Conference (ACC)},
  pages={1501--1506},
  year={2019},
  organization={IEEE}
}

@book{rockafellar1998variational,
  title={Variational analysis},
  author={Rockafellar, R Tyrrell and Wets, Roger JB},
  year={1998},
  publisher={Springer}
}

@article{le2018dc,
  title={DC programming and DCA: thirty years of developments},
  author={Le Thi, Hoai An and Pham Dinh, Tao},
  journal={Mathematical Programming},
  volume={169},
  number={1},
  pages={5--68},
  year={2018},
  publisher={Springer}
}

@article{rockafellar2002conditional,
  title={Conditional value-at-risk for general loss distributions},
  author={Rockafellar, R Tyrrell and Uryasev, Stanislav},
  journal={Journal of banking \& finance},
  volume={26},
  number={7},
  pages={1443--1471},
  year={2002},
  publisher={Elsevier}
}

@article{acerbi2002coherence,
  title={On the coherence of expected shortfall},
  author={Acerbi, Carlo and Tasche, Dirk},
  journal={Journal of Banking \& Finance},
  volume={26},
  number={7},
  pages={1487--1503},
  year={2002},
  publisher={Elsevier}
}

@book{rockafellar1970convex,
  title={Convex Analysis},
  author={Rockafellar, R Tyrrell},
  year={1970},
  publisher={Princeton University Press},
  address={Princeton, NJ}
}
\newpage
\begin{APPENDICES}
\section{Examples}
\label{appendix_example}

\begin{example}\label{ex:DC_procedure}
Consider a DRCCP~\eqref{eq_drccp_org_ch2} whose reference distribution \(\hat{\Pr}\) is the empirical distribution of $N=3$ scenarios and $p=2$. Let the risk level be $\varepsilon=0.5$, Wasserstein radius $\theta=0.2$, cost vector $\bm{c} = (-1,-1)^\top$, deterministic feasible region $\X = [0,10]^2$. The constraint functions are given by $\bm a_i(\x)=\x, b_i(\x)=b =1, I=1$ with scenarios $\hat{\rxi}^1=(1,0)^\top$, $\hat{\rxi}^2=(0,1)^\top$, and $\hat{\rxi}^3=(2,3)^\top$.
Applying SCA Algorithm~\ref{alg_DCA} with an initial solution $\bm{x}^{(0)} = \bm 0$ and a tolerance $\tau = 10^{-3}$, we obtain:
\begin{align*}
v^{(0)}(\mathrm{DC}) = -0.485, v^{(1)}(\mathrm{DC}) = -0.742,v^{(2)}(\mathrm{DC}) = -0.742. 
\end{align*}
Since $|v^{(2)}(\mathrm{DC})-v^{(1)}(\mathrm{DC})|<\tau$, the algorithm terminates after 3 iterations and yields an objective value $v^{\mathrm{DC}} = -0.742$ with a solution $\bm{x}^{\mathrm{DC}} = (0.371, 0.371)$. 

In this example, if we replace the right-hand CVaR function with $0$, we obtain CVaR approximation~\eqref{worse_case_cvar_q_formulation_a}, which
admits the following form:
\begin{align*}
v^{\CVaR} = \min_{\bm{x}\in[0,10]^2,\beta\in\mathbb{R}}\left\{ \bm{c}^{\top}\bm{x} : \theta\varepsilon^{-1}\|\bm{x}\|_2 +  \beta +\frac{1}{N\varepsilon}\sum_{j=1}^{N} \left( \hat{\rxi}^{j\top}\bm{x} - 1 -\beta \right)_+ \leq 0 \right\}
\end{align*}
with an optimal objective value $v^{\CVaR} = -0.485$ and an optimal solution $\bm x^{\CVaR}=(0.419, 0.066)$.  Clearly, in this example, we have $v^{\mathrm{DC}} < v^{\CVaR}$.
\QEDB
\end{example}

\begin{example}
\label{example_cvar_infeasible}
Consider a DRCCP~\eqref{eq_drccp_org_ch2} whose reference distribution \(\hat{\Pr}\) is the empirical distribution of $N=3$ scenarios and $p=2$. Let the risk level be $\varepsilon=0.5$, Wasserstein radius $\theta=0.1$, cost vector $\bm{c} = (-1,-1)^\top$, deterministic feasible region $\X = [1,2]^2$. The constraint functions are given by $\bm a_i(\x)=\x, b_i(\x)=b=3, I=1$ with scenarios $\hat{\rxi}^1=(2,0)^\top$, $\hat{\rxi}^2=(0,2)^\top$, and $\hat{\rxi}^3=(5/2,5/2)^\top$.
Let $\hat{\bm x}\in\argmin_{\bm x\in\X}\hat{\Pr}\text{-}\CVaR_{1-\varepsilon}[ \{\bm x^\top\hat{\rxi}-b\}_+]$. A straightforward calculation shows that $\hat{\bm x}=(1,1)$. Applying SCA Algorithm~\ref{alg_DCA} with an initial solution $\bm{x}^{(0)} = (1,1)$ and a tolerance $\tau = 10^{-3}$, we obtain $v^{\mathrm{DC}} = -2.106$ with a solution $\bm x^{\mathrm{DC}} = (1.053, 1.053)$. The corresponding CVaR approximation~\eqref{worse_case_cvar_q_formulation_a} can be written as
\begin{align*}
v^{\CVaR} = \min_{\bm{x}\in [1,2]^2,\beta\in\mathbb{R}}\left\{ \bm{c}^{\top}\bm{x} : \theta\varepsilon^{-1}\|\bm{x}\|_2 +  \beta +\frac{1}{N\varepsilon}\sum_{j=1}^{N} \left( \hat{\rxi}^{j\top}\bm{x} - 3 -\beta\right)_+ \leq 0 \right\}.
\end{align*}
However, this CVaR approximation is infeasible over the entire set $\bm x\in[1,2]^2$. 
Therefore, in this example, SCA Algorithm~\ref{alg_DCA} finds a feasible solution to the corresponding DRCCP, while CVaR approximation~\eqref{worse_case_cvar_q_formulation_a} fails to find any solution.
\QEDB
\end{example}

\begin{example}
\label{example_revisit_SCA}
Let us revisit Example~\ref{ex:DC_procedure} with Wasserstein radius $\theta=0.5$, while keeping all other parameters unchanged. Applying SCA Algorithm~\ref{alg_DCA} with an initial solution $\bm{x}^{(0)} = \bm{0}$ and a tolerance $\tau=10^{-3}$, we obtain:
\begin{align*}
v^{(0)}(\mathrm{DC}) = -0.397, \quad v^{(1)}(\mathrm{DC}) = -0.397.
\end{align*}
Since $|v^{(1)}(\mathrm{DC})-v^{(0)}(\mathrm{DC})|<\tau$, the algorithm converges in 2 iterations and provides an objective value $v^{\mathrm{DC}} = -0.397$ with a solution $\bm{x}^{\mathrm{DC}} = (0.247, 0.150)$. Correspondingly, CVaR approximation~\eqref{worse_case_cvar_q_formulation_a} yields $v^{\CVaR} = -0.397$ with $\bm{x}^{\CVaR} = (0.247, 0.150)$. Clearly, $v^{\mathrm{DC}} = v^{\CVaR}$ and $\bm{x}^{\mathrm{DC}} =\bm{x}^{\CVaR}$.
This result occurs because $\bm{x}^{\CVaR}$ minimizes $\hat{\Pr}\text{-}\CVaR_{1-\varepsilon}[ \max_{i\in[I]} \{\bm{a}_i(\bm{x})^\top\hat{\rxi}-b_i(\bm{x})\}_+ ]$, i.e., $\bm{x}^{\CVaR} \in \argmin_{\bm{x}\in\X} \hat{\Pr}\text{-}\CVaR_{1-\varepsilon}[ \max_{i\in[I]} \{\bm{a}_i(\bm{x})^\top\hat{\rxi}-b_i(\bm{x})\}_+ ]$.
Therefore, SCA Algorithm~\ref{alg_DCA} converges to the same solution as CVaR approximation~\eqref{worse_case_cvar_q_formulation_a} after the first iteration.
\QEDB
\end{example}

\begin{example}
\label{example_compare_also_a}
Let us revisit Example~\ref{ex:DC_procedure} to compare the performance of SCA Algorithm~\ref{alg_DCA} and ALSO-X\# Algorithm \ref{alg_alsox_sharp}. In this example, SCA Algorithm~\ref{alg_DCA} converges in $3$ iterations and yields an objective value $v^{\mathrm{DC}} = -0.742$ with the solution $\bm{x}^{\mathrm{DC}} = (0.371, 0.371)$. Recall that CVaR approximation yields $v^{\CVaR} = -0.485$ with the solution $\bm{x}^{\CVaR} = (0.419, 0.066)$. We now apply ALSO-X\# Algorithm~\ref{alg_alsox_sharp} to the same instance. Implemented with tolerance $\delta_1 = 10^{-3}$ and initial bounds $t_L = -5$, $t_U = 0$, the algorithm converges to an objective value $v^{\as} = -0.522$ with the solution $\bm{x}^{\as} = (0.451, 0.071)$ after $13$ iterations. We obtain
\[
v^{\mathrm{DC}} = -0.742 < v^{\as} = -0.522 < v^{\CVaR} = -0.485.
\]
This confirms that SCA Algorithm \ref{alg_DCA} can outperform both ALSO-X\# Algorithm  \ref{alg_alsox_sharp} and CVaR approximation~\eqref{worse_case_cvar_q_formulation_a}.
\QEDB
\end{example}

\begin{example}
\label{example_compare_also_b}
Consider a DRCCP~\eqref{eq_drccp_org_ch2} whose reference distribution \(\hat{\Pr}\) is the empirical distribution of $N=5$ scenarios and $p=2$. Let the risk level be $\varepsilon=0.5$, Wasserstein radius $\theta=0.1$, cost vector $\bm{c} = (-2,-2)^\top$, deterministic feasible region $\X = [0,10]^2$. The constraint functions are given by $\bm a_i(\x)=\x, b_i(\x)=b=3, I=1$ with scenarios $\hat{\rxi}^1=(1,3)^\top$, $\hat{\rxi}^2=(1,9)^\top$, $\hat{\rxi}^3=(7,9)^\top$, $\hat{\rxi}^4=(10,4)^\top$, and $\hat{\rxi}^5=(8,2)^\top$.
Applying SCA Algorithm~\ref{alg_DCA} with an initial solution $\bm{x}^{(0)} = \bm{0}$ and a tolerance $\tau = 10^{-3}$, we obtain:
\begin{align*}
v^{(0)}(\mathrm{DC}) = -0.840, v^{(1)}(\mathrm{DC}) = -0.939, v^{(2)}(\mathrm{DC}) = -0.939.
\end{align*}
Since $|v^{(2)}(\mathrm{DC})-v^{(1)}(\mathrm{DC})|<\tau$, the algorithm converges in 3 iterations with the objective value $v^{\mathrm{DC}} = -0.939$ and the solution $\bm{x}^{\mathrm{DC}} = (0.168, 0.302)$. The corresponding CVaR approximation~\eqref{worse_case_cvar_q_formulation_a} yields:
\begin{align*}
v^{\CVaR} = \min_{\bm{x}\in[0,10]^2,\beta\in\mathbb{R}}\left\{ \bm{c}^{\top}\bm{x} : \theta\varepsilon^{-1}\|\bm{x}\|_2 +  \beta+\frac{1}{N\varepsilon}\sum_{j=1}^{N} \left( \hat{\rxi}^{j\top}\bm{x} - 3 -\beta\right)_+ \leq 0 \right\}
\end{align*}
with an optimal value $v^{\CVaR} = -0.840$ and the solution $\bm{x}^{\CVaR} = (0.209, 0.211)$. Now applying ALSO-X\# Algorithm~\ref{alg_alsox_sharp} with a tolerance $\delta_1 = 10^{-3}$ and initial bounds $t_L = -5$, $t_U = 0$, the algorithm converges in $13$ iterations to an objective value $v^{\as} = -1.049$ with the solution $\bm{x}^{\as} = (0.262, 0.263)$. We observe that
\[
v^{\as} = -1.049 < v^{\mathrm{DC}} = -0.939 < v^{\CVaR} = -0.840.
\]
This demonstrates that ALSO-X\# Algorithm~\ref{alg_alsox_sharp} can outperform both SCA Algorithm~\ref{alg_DCA} and CVaR approximation~\eqref{worse_case_cvar_q_formulation_a}, and hence that SCA Algorithm~\ref{alg_DCA} and ALSO-X\# Algorithm~\ref{alg_alsox_sharp} are generally not comparable.
\QEDB
\end{example}

\begin{example}
\label{example_sensitivity_detail}
We revisit Example~\ref{example_sensitivity} and analyze the performance of SCA Algorithm~\ref{alg_DCA}, CVaR approximation~\eqref{worse_case_cvar_q_formulation_a},
and ALSO-X\# Algorithm~\ref{alg_alsox_sharp} as the parameter
$\alpha\in[0.7,0.9]$ varies in the empirical distribution. Recall that we set $\varepsilon=0.5$, $\theta=0.01$, $p=2$,  $g(\bm x,\hat{\rxi})=\bm x^\top\hat{\rxi}+1$ and $N=5$. 
For the corresponding CVaR approximation~\eqref{worse_case_cvar_q_formulation_a}, the optimal solution admits a closed-form characterization. Under the five equally weighted empirical scenarios, the corresponding losses are $1-x_1$, $1-\alpha x_2$, $1-x_3$, $1-x_1-x_2-x_3$, and $1-x_1-x_2-x_3$. Since the last two losses are no larger than the first three on $\X$, the constraint of CVaR approximation~\eqref{worse_case_cvar_q_formulation_a} reduces to
\begin{equation}
    0.02\|\bm x\|_2+1-0.4(x_1+\alpha x_2+x_3)+0.2\max\{x_1,\alpha x_2,x_3\}\le 0,
\label{exp7_cvar_constraint}
\end{equation}
and an optimal solution can be chosen in the form
\[
\bm x=u(1,r,1)^\top,\quad u\ge 0, r\ge 0.
\]
For $\alpha\in[0.7,0.9]$, no ratio $r>1/\alpha$ improves the objective value $\bm c^\top\bm x$. For $0\le r\le 1/\alpha$, by constraint~\eqref{exp7_cvar_constraint}, we know that the smallest feasible value of $u$ is
$u(r,\alpha)=1/(0.6+0.4\alpha r-0.02\sqrt{2+r^2})$. Consequently, we have
\begin{equation}
    v^{\CVaR}(\alpha)=\min_{0\le r\le 1/\alpha}\frac{2+r}{0.6+0.4\alpha r-0.02\sqrt{2+r^2}}.
\label{exp7_cvar_alpha}
\end{equation}
To characterize the optimal ratio $r_\alpha$ on $[0,1/\alpha]$ in equation~\eqref{exp7_cvar_alpha}, we define $\alpha_L=(0.6-0.02\sqrt{2})/0.8 \approx 0.715$, and let $\alpha_U\in(0.75,0.8)$ denote the unique solution of $0.6-0.8\alpha_U+0.04(1-\alpha_U)/\sqrt{2\alpha_U^2+1}=0$, which is $\alpha_U\approx 0.758$. Further, let
$q_\alpha=20\alpha-15$. Then we have
\[
r_\alpha=
\begin{cases}
0, & 0.7\leq\alpha\leq\alpha_L,\\[3pt]
\displaystyle\frac{1+q_\alpha\sqrt{3-2q_\alpha^2}}{1-q_\alpha^2},
& \alpha_L<\alpha<\alpha_U,\\[7pt]
1/\alpha, & \alpha_U\leq\alpha\leq0.9,
\end{cases}
\qquad
u_\alpha=\frac{1}{0.6+0.4\alpha r_\alpha-0.02\sqrt{2+r_\alpha^2}}.
\]
The optimal solution and objective value of CVaR approximation~\eqref{worse_case_cvar_q_formulation_a} are therefore
\[
\bm x^{\CVaR}(\alpha)=u_\alpha(1,r_\alpha,1)^\top,
\qquad
v^{\CVaR}(\alpha)=(2+r_\alpha)u_\alpha.
\]
Thus, the solution varies continuously with $\alpha$.
In particular, at $\alpha=0.7$, we obtain  $v^{\CVaR}(0.7)=3.498$ with
$\bm x^{\CVaR}(0.7)=(1.749,0.000,1.749)^\top$, at $\alpha=0.74$, we obtain  $v^{\CVaR}(0.74)=3.481$ with $\bm x^{\CVaR}(0.74)=(1.296,0.889,1.296)^\top$, while at $\alpha=0.9$, we have $v^{\CVaR}(0.9)=3.227$ with $\bm x^{\CVaR}(0.9)=(1.037,1.153,1.037)^\top$. These solutions agree exactly with the numerical results.

ALSO-X\# Algorithm~\ref{alg_alsox_sharp} solves a sequence of lower-level problems and enforces feasibility through the checking condition \eqref{alsoxs_drccp_formualtion_upper}. In this instance, the numerical outputs vary with $\alpha$ and can be directly compared with solutions of the CVaR approximation discussed above. 
Specifically, for $\alpha=0.7$, it returns $\bm x^{\as}=(1.049,0.000,1.049)^\top$ with
$v^{\as}=2.099$, whereas for $\alpha=0.9$, it returns $\bm x^{\as}=(1.038,1.153,1.038)^\top$ with $v^{\as}=3.228$.
These results show that the solutions returned by ALSO-X\# vary with $\alpha$.

In contrast, SCA Algorithm~\ref{alg_DCA} returns the same solution for all $\alpha\in[0.7,0.9]$. Initializing the algorithm at $\bm{x}^{(0)}=(10000,0,0)^\top$, we find that only the second and third losses are positive at this point, so we have
\[
g_2(\bm{x}^{(0)})=0.8, \qquad
\bm{s}_2^{(0)}=(0,-0.4\alpha,-0.4)^\top.
\]
Therefore, the first subproblem~\eqref{eq_dc} has the constraint
\[
0.02\|\bm{x}\|_2+0.2-0.4x_1+0.2\max\{x_1,\alpha x_2,x_3\}\leq0.
\]
Since $\max\{x_1,\alpha x_2,x_3\}\ge x_1$, every point feasible to this constraint satisfies $x_1\geq1+0.1\|\bm{x}\|_2\geq1+0.1x_1$. Then the objective value satisfies $\bm c^\top\bm x =x_1+x_2+x_3\geq x_1\geq10/9$. The point $\bm{x}^{(1)}=(10/9,0,0)^\top$ attains this bound and is feasible to the first subproblem. It is therefore its unique optimal solution. At $\bm{x}^{(1)}$, the next subproblem~\eqref{eq_dc} is identical and returns $\bm{x}^{(2)}=\bm{x}^{(1)}$. Thus, with tolerance $\tau=10^{-3}$, SCA Algorithm~\ref{alg_DCA} converges in two iterations to $\bm x^{\mathrm{DC}}=({10}/{9},0,0)^\top, v^{\mathrm{DC}}={10}/{9}$. Notably, neither $\bm x^{\mathrm{DC}}$ nor $v^{\mathrm{DC}}$ depends on $\alpha$.

In this setting, the scenario $\hat{\rxi}^2$ falls in Region III of Figure~\ref{fig:tail_decomposition}, which corresponds to the term $\hat{\Pr}\text{-}\CVaR_{1-\varepsilon}[\max_{i\in[I]}\{\bm{a}_i(\bm{x})^\top\hat{\rxi}-b_i(\bm{x})\}_+]$ in DRCCP~\eqref{eq_dc_orginial}. The key insight is that DC reformulation~\eqref{eq_dc_orginial} explicitly separates and offsets the positive-tail contribution arising from Region III through the term $\hat\Pr\text{-}\CVaR_{1-\varepsilon}[(\hat{\rxi}^{\top}\bm{x} + 1)_+]$, and SCA subproblem~\eqref{eq_dc} uses this structure via its linearization,
making it closer to the original DRCCP formulation~\eqref{1_wasser_reformualtion}. This exclusion ensures that scenarios in Region III do not adversely affect the solution. Consequently, when $\alpha$ changes from $0.7$ to $0.9$, SCA Algorithm~\ref{alg_DCA} remains unaffected. However, both CVaR approximation~\eqref{worse_case_cvar_q_formulation_a} and the lower-level problem of ALSO-X\# Algorithm~\ref{alg_alsox_sharp} retain the positive-tail contribution (i.e., $\hat\Pr\text{-}\CVaR_{1-\varepsilon}[(\hat{\rxi}^{\top}\bm{x} + 1)_+]$ in DRCCP~\eqref{eq_dc_orginial}), which can make their candidate solutions sensitive to changes in Region III. As a result, when the scenario becomes more extreme ($\alpha$ increases), both methods remain sensitive to changes in Region III; in particular, ALSO-X\# deteriorates as $\alpha$ increases, while CVaR approximation remains substantially conservative over the entire range. This behavior illustrates a limitation of existing convex approximations when feasibility is driven by a subset of scenarios (here, as $\alpha$ increases from $0.7$ to $0.9$). The proposed SCA framework, however, remains unchanged under the same perturbation.
\QEDB
\end{example}

\section{Convergence and Iteration Complexity of SCA Algorithm}
\label{appendix_complexity}

\subsection{Convergence Analysis}
This subsection proves Proposition~\ref{prop_converge} in three steps. We first establish two basic properties of the iterates of SCA Algorithm~\ref{alg_DCA}, which give the first two claims of the proposition. We then show that Assumption~\ref{assumption_lambda} guarantees Slater's condition for every SCA subproblem~\eqref{eq_dc} (Lemma~\ref{lemma_dc}). Finally, we combine these results to show that every cluster point of the iterates is a critical point of DRCCP~\eqref{eq_dc_orginial}.

We begin with the two basic properties, neither of which requires Assumption~\ref{assumption_lambda}:
\begin{enumerate*}[label=(\roman*)]
\item\label{property_feasible} all iterates $\bm x^{(k)}$ remain feasible, i.e., $\bm x^{(k)} \in \Omega(\mathrm{DC})$ for all $k\in \Ze_+$, because the linearization of the convex function $g_2(\cdot)$ at $\bm x^{(k)}$ is a global lower bound on $g_2(\cdot)$ and hence $\Omega(\mathrm{DC})_{\bm x^{(k)}}\subseteq\Omega(\mathrm{DC})$; and
\item\label{property_monotone} the sequence of objective values $\{\bm c^\top \bm x^{(k)}\}_{k\in \Ze_+}$ is monotonically nonincreasing and bounded below, and is therefore convergent, because $\bm x^{(k)}$ is feasible for the $k$-th subproblem~\eqref{eq_dc} and $\X$ is compact.
\end{enumerate*}

To show that every cluster point of the iterates is a critical point of DRCCP~\eqref{eq_dc_orginial}, we additionally impose Assumption~\ref{assumption_lambda} stated in Section~\ref{section_sub_sca}, which can equivalently be written in multiplier form: for every $\bm x\in\Omega(\mathrm{DC})$, if $\lambda\geq 0$ satisfies $0\in\lambda(\partial g_1(\bm x)-\partial g_2(\bm x))+\N_\X(\bm x)$ and $\lambda[g_1(\bm x)-g_2(\bm x)]=0$, then $\lambda=0$. Indeed, $\lambda>0$ forces the constraint to be active, and since $\N_\X(\bm x)$ is a cone, dividing by $\lambda$ yields $\bm 0\in\partial g_1(\bm x)-\partial g_2(\bm x)+\N_\X(\bm x)$. In particular, $\lambda$ can be viewed as the Lagrange multiplier associated with the constraint $g_1(\bm x)-g_2(\bm x)\le 0$, so Assumption~\ref{assumption_lambda} rules out nonzero abnormal multipliers. As shown in Lemma~\ref{lemma_dc} below, it ensures the relative Slater's condition for the convexified feasible set generated at any feasible point. When $g_1(\cdot)$ and $g_2(\cdot)$ are differentiable, Assumption~\ref{assumption_lambda} states that $-(\nabla g_1(\bm x)-\nabla g_2(\bm x))\notin\N_\X(\bm x)$ at every active feasible point $\bm x$, or equivalently, since $\N_\X(\bm x)$ is the polar cone of the tangent cone of $\X$ at $\bm x$, that some tangent direction $\bm d$ satisfies $(\nabla g_1(\bm x)-\nabla g_2(\bm x))^\top\bm d<0$. By Theorem~6.14 and Exercise~6.39 of \citet{rockafellar1998variational}, this is the normal-cone constraint qualification for the system $\bm x\in\X$, $g_1(\bm x)-g_2(\bm x)\le 0$, which is equivalent to its tangent-cone form because $\X$ is convex; when $\X=\Re^n$, it is the Mangasarian--Fromovitz constraint qualification (MFCQ; Example~6.40 therein). \citet{hong2011sequential} impose the same type of condition on the chance constraint itself (their Assumption~6) to derive Slater's condition for their SCA subproblems (their Lemma~4), and Lemma~\ref{lemma_dc} below follows a similar argument for the DC constraint. Assumption~\ref{assumption_lambda} can therefore be interpreted as a nonsmooth counterpart of these classical constraint qualifications, with the gradient of $g_1(\cdot)-g_2(\cdot)$ replaced by the generalized difference $\partial g_1(\bm x)-\partial g_2(\bm x)$, following the subdifferential framework of nonsmooth DC programming (see, e.g., \citealt{tao1997convex, le2018dc}). %

To establish convergence to critical points, we need Slater's condition to hold for the convex SCA subproblem~\eqref{eq_dc} generated at each feasible point, which we derive from Assumption~\ref{assumption_lambda} in Lemma~\ref{lemma_dc} below. Since the convex set $\X$ may have an empty interior, we state Slater's condition in terms of its relative interior \(\operatorname{ri}(\mathcal{X})\), which is always nonempty. For any
\(\bm y\in \Omega(\mathrm{DC})\) and \(\bm s_2(\bm y)\in \partial g_2(\bm y)\), recall that
\(
\Omega(\mathrm{DC})_{\bm y} =
\left\{
\bm x\in \X:
g_1(\bm x)-\left[g_2(\bm y)+\bm s_2(\bm y)^\top(\bm x-\bm y)\right]\le 0
\right\},
\)
for which Slater's condition (see, e.g., \citealt{slater1950lagrange}) requires the existence of
some \(\bm{\tilde x}\in \operatorname{ri}(\X)\) such that
\begin{equation}\label{eq_slater}
g_1(\bm{\tilde x})-\left[g_2(\bm y)+\bm s_2(\bm y)^\top(\bm{\tilde x}-\bm y)\right] <0.
\end{equation}

\begin{lemma}
\label{lemma_dc}
Under Assumption~\ref{assumption_lambda}, 
for every \(\bm y\in \Omega(\mathrm{DC})\) and every
\(\bm s_2(\bm y)\in \partial g_2(\bm y)\), the convexified feasible set
\(\Omega(\mathrm{DC})_{\bm y}\) satisfies Slater's condition relative to \(\X\); that is,
there exists \(\bm{\tilde x}\in \operatorname{ri}(\X)\) satisfying~\eqref{eq_slater}.
\end{lemma}
\begin{proof}
Suppose by contradiction that there exist \(\bm y\in \Omega(\mathrm{DC})\)
and \(\bm s_2(\bm y)\in \partial g_2(\bm y)\) such that Slater's condition fails for
\(\Omega(\mathrm{DC})_{\bm y}\). Then, for every \(\bm x\in \operatorname{ri}(\X)\), we have
$g_1(\bm x)-[g_2(\bm y)+\bm s_2(\bm y)^\top(\bm x-\bm y)]\ge 0$.
Since \(\X\) is convex with a nonempty relative interior, the closure of
\(\operatorname{ri}(\X)\) is \(\X\). Moreover, \(g_1(\cdot)\) is continuous on \(\X\), and
\(g_2(\bm y)+\bm s_2(\bm y)^\top(\bm x-\bm y)\) is affine in \(\bm x\). Therefore, by continuity, this inequality holds for every \(\bm x\in \X\). On the other hand, since \(\bm y\in \Omega(\mathrm{DC})\), we have $g_1(\bm y)-g_2(\bm y)\le 0$. Combining this inequality with the previous one evaluated at \(\bm x=\bm y\), we obtain $g_1(\bm y)-g_2(\bm y)=0$. Furthermore, by $g_1(\bm x)-[g_2(\bm y)+\bm s_2(\bm y)^\top(\bm x-\bm y)]$ being nonnegative on $\X$ and vanishing at $\bm x=\bm y$,
 $\bm y\in \argmin_{\bm x\in \X} \{ g_1(\bm x)-[g_2(\bm y)+\bm s_2(\bm y)^\top(\bm x-\bm y)] \}$.
Based on the necessary optimality condition for this convex optimization problem, there exists \(\bm s_1(\bm y)\in \partial g_1(\bm y)\) such that
$-[\bm s_1(\bm y)-\bm s_2(\bm y)]\in \N_\X(\bm y)$.
Since $\bm s_2(\bm y)\in\partial g_2(\bm y)$, we have $\bm 0\in\partial g_1(\bm y)-\partial g_2(\bm y)+\N_\X(\bm y)$ at a point $\bm y\in\Omega(\mathrm{DC})$ with $g_1(\bm y)=g_2(\bm y)$, contradicting Assumption~\ref{assumption_lambda}. Therefore, Slater's condition
must hold for \(\Omega(\mathrm{DC})_{\bm y}\). This completes the proof.
\QEDA
\end{proof}
Based on Lemma~\ref{lemma_dc} and the discussions above, we now prove Proposition~\ref{prop_converge} stated in Section~\ref{section_sub_sca}.

\begin{proof}
  The first two claims of the proposition are properties~\ref{property_feasible} and~\ref{property_monotone} above. To prove the last claim, we borrow ideas from the proof of Property~3 in \citet{hong2011sequential}, who show for continuously differentiable $g_1$ and $g_2$ that every cluster point $\bar{\bm x}$ of the iterates is optimal for the subproblem linearized at $\bar{\bm x}$ and is therefore a KKT point under Slater's condition. The main difficulty in our case is that $g_2$ is~\emph{nonsmooth} and each subproblem~\eqref{eq_dc} uses an arbitrary subgradient, so the linearizations of $g_2$ need not converge along a subsequence of iterates that converges to $\bar{\bm x}$.

  Let $\bar{\bm x}$ be a cluster point of $\{\bm x^{(k)}\}_{k\in\Ze_+}$ and let $v_\infty:=\lim_{k\to\infty}\bm c^\top\bm x^{(k)}$, so that $\bar{\bm x}\in\Omega(\mathrm{DC})$ and $\bm c^\top\bar{\bm x}=v_\infty$ by continuity. We proceed in three steps: Step~1 obtains a subgradient $\bar{\bm s}_2\in\partial g_2(\bar{\bm x})$ as a limit of the subgradients $\bm s_2^{(k)}$ that define the subproblems~\eqref{eq_dc}, Step~2 shows that $\bar{\bm x}$ is optimal for the subproblem linearized at $\bar{\bm x}$ with $\bar{\bm s}_2$, and Step~3 concludes that $\bar{\bm x}$ is a critical point of DRCCP~\eqref{eq_dc_orginial}.

  \begin{enumerate}[label=Step~\arabic*., leftmargin=*]
  \item Since the subgradients of the finite convex function $g_2$ are bounded on $\X$ \citep[\,Theorem~24.7]{rockafellar1970convex}, we can select a subsequence such that both $\bm x^{(k_l)}\to\bar{\bm x}$ and $\bm s_2^{(k_l)}\to\bar{\bm s}_2$ for some $\bar{\bm s}_2$, and taking the limit in the subgradient inequality $g_2(\bm x)\geq g_2(\bm x^{(k_l)})+(\bm s_2^{(k_l)})^\top(\bm x-\bm x^{(k_l)})$ shows that $\bar{\bm s}_2\in\partial g_2(\bar{\bm x})$.
  \item We show that $\bar{\bm x}$ is optimal for the subproblem linearized at $\bar{\bm x}$ with the subgradient $\bar{\bm s}_2$, for which it is feasible because $g_1(\bar{\bm x})\leq g_2(\bar{\bm x})$. Suppose by contradiction that some $\bm x'$ feasible for this subproblem satisfies $\bm c^\top\bm x'<v_\infty$. By Lemma~\ref{lemma_dc} applied at $\bm y=\bar{\bm x}$ with $\bm s_2(\bm y)=\bar{\bm s}_2$, some $\tilde{\bm x}\in\X$ satisfies the constraint of this subproblem strictly, and so does $\bm x_t:=(1-t)\bm x'+t\tilde{\bm x}$ for every $t\in(0,1]$ because $g_1$ is convex; we fix $t$ small enough such that $\bm c^\top\bm x_t<v_\infty$. Since the constraint of the $k_l$-th subproblem~\eqref{eq_dc} converges to this constraint as $l\to\infty$, $\bm x_t$ is feasible for the $k_l$-th subproblem for all sufficiently large $l$, so $\bm c^\top\bm x^{(k_l+1)}\leq\bm c^\top\bm x_t<v_\infty$, which contradicts property~\ref{property_monotone}. Hence, $\bar{\bm x}$ is optimal for this subproblem.
  \item Since the convex subproblem in Step~2 satisfies Slater's condition by Lemma~\ref{lemma_dc}, its KKT conditions hold at its optimal solution $\bar{\bm x}$ \citep[\,Theorems~28.2 and~28.3]{rockafellar1970convex}, i.e., there exist $\lambda\geq 0$ and $\bm s_1\in\partial g_1(\bar{\bm x})$ such that
  \[
  \bm 0\in\bm c+\lambda(\bm s_1-\bar{\bm s}_2)+\N_\X(\bar{\bm x}),\qquad \lambda[g_1(\bar{\bm x})-g_2(\bar{\bm x})]=0,
  \]
  which matches the definition of a critical point of DRCCP~\eqref{eq_dc_orginial} with $\bm s_2=\bar{\bm s}_2$.
  \end{enumerate}
\QEDA
\end{proof}

\subsection{Iteration Complexity for Approximate Criticality}
After establishing the convergence of SCA Algorithm~\ref{alg_DCA}, we now study its iteration complexity for solving DRCCP~\eqref{eq_dc_orginial}. In particular, we quantify the number of iterations required to obtain an approximate critical point. Recall that in SCA Algorithm~\ref{alg_DCA}, $v^{(k)}(\mathrm{DC})$ denotes the optimal objective value of the $k$-th subproblem \eqref{eq_dc}. By construction of SCA Algorithm~\ref{alg_DCA}, we have $v^{(k)}(\mathrm{DC})=\bm c^\top \bm x^{(k+1)}$. For $k\ge 1$, we use $\Delta_k\ge 0$ to denote the objective improvement attained by the $k$-th iteration: 
\begin{subequations}
\begin{equation}\label{eq:Delta_k}
\Delta_k := v^{(k-1)}(\mathrm{DC})-v^{(k)}(\mathrm{DC})
= \bm c^\top\bm x^{(k)}-\bm c^\top\bm x^{(k+1)}.
\end{equation}
Since $\bm x^{(k+1)}$ is an optimal solution of the $k$-th subproblem \eqref{eq_dc}, it follows that for any $\bm x\in \Omega(\mathrm{DC})_{\bm x^{(k)}}$, we have $\bm c^\top\bm x \ge v^{(k)}(\mathrm{DC})$, which implies that
\begin{equation}\label{eq:approx_opt_property}
\bm c^\top\bm x - \bm c^\top\bm x^{(k)} \ge -\Delta_k,\quad \forall\, \bm x\in \Omega(\mathrm{DC})_{\bm x^{(k)}}.
\end{equation}
\end{subequations}
For any feasible point $\bar{\bm x}\in\Omega(\mathrm{DC})$, we define the criticality gap at \(\bar{\bm x}\) as
\[
\mathcal G(\bar{\bm x}):=\bm c^\top\bar{\bm x}-\sup_{\bm s_2(\bar{\bm x})\in\partial g_2(\bar{\bm x})}\min_{\bm x\in\Omega(\mathrm{DC})_{\bar{\bm x}}}\bm c^\top\bm x,
\]
where the supremum is over the subgradient $\bm s_2(\bar{\bm x})$ that defines $\Omega(\mathrm{DC})_{\bar{\bm x}}$, so that $\mathcal G(\bar{\bm x})$ does not depend on a particular choice of subgradient, and $\mathcal G(\bar{\bm x})\geq 0$ because $\bar{\bm x}\in\Omega(\mathrm{DC})_{\bar{\bm x}}$.
A feasible point $\bar{\bm x}$ is called an $\bar{\epsilon}$-approximate critical point measured by the criticality gap if $\mathcal G(\bar{\bm x})\le \bar{\epsilon}$. In particular, under Assumption~\ref{assumption_lambda}, $\mathcal G(\bar{\bm x})=0$ if and only if $\bar{\bm x}$ is a critical point of DRCCP~\eqref{eq_dc_orginial}. Indeed, if $\bar{\bm x}$ is critical with some $\bm s_2\in\partial g_2(\bar{\bm x})$, then it satisfies the KKT conditions of the convex subproblem linearized at $\bar{\bm x}$ with $\bm s_2$ and is thus optimal for it, so $\mathcal G(\bar{\bm x})=0$. Conversely, if $\mathcal G(\bar{\bm x})=0$, then repeating Step~2 of the proof of Proposition~\ref{prop_converge} along a maximizing sequence of subgradients in the compact set $\partial g_2(\bar{\bm x})$ shows that $\bar{\bm x}$ is optimal for the subproblem linearized at $\bar{\bm x}$ with a limit of this sequence, and Step~3 shows that $\bar{\bm x}$ is critical.
More generally, a feasible point with a small criticality gap is referred to as approximately critical for the original DC problem.
Based on the definitions above, we now prove Theorem~\ref{thm:iteration-complexity} stated in Section~\ref{section_sub_sca}.

\begin{proof}
By equation \eqref{eq:Delta_k} and telescoping, we have
\[
\sum_{i=1}^{k}\Delta_i
=\sum_{i=1}^{k}\bigl(v^{(i-1)}(\mathrm{DC})-v^{(i)}(\mathrm{DC})\bigr)
=v^{(0)}(\mathrm{DC})-v^{(k)}(\mathrm{DC})
=\bm c^\top\bm x^{(1)}-\bm c^\top\bm x^{(k+1)}
\le \bar M,
\]
which implies $\min_{1\le i\le k}\Delta_i \le {\bar M}/{k}$. Moreover, since $\bm s_2^{(i)}$ is one of the subgradients over which the supremum in $\mathcal G(\bm x^{(i)})$ is taken, \eqref{eq:approx_opt_property} implies $\mathcal G(\bm x^{(i)})\leq\Delta_i$ for each $i$, and hence $\min_{1\le i\le k}\mathcal G(\bm x^{(i)})\leq\min_{1\le i\le k}\Delta_i\leq\bar M/k$.
\QEDA
\end{proof}
Theorem~\ref{thm:iteration-complexity} implies that after $k+1$ iterations of SCA Algorithm~\ref{alg_DCA}, we can identify an iterate $\bm{x}^{(\ell)}$ whose criticality gap satisfies $\mathcal G(\bm x^{(\ell)})\le \bar{\epsilon}$, with $\bar{\epsilon} = \O(1/k)$. Since, under Assumption~\ref{assumption_lambda}, a zero criticality gap implies criticality of DRCCP~\eqref{eq_dc_orginial} under the relative Slater condition, this result provides a finite iteration guarantee for approaching criticality of DRCCP~\eqref{eq_dc_orginial} as measured by $\mathcal G (\cdot)$.
It also provides a theoretical guarantee on the convergence rate of our algorithm: To obtain an $\bar{\epsilon}$-approximate critical point, which is guaranteed once $\Delta_\ell\le\bar{\epsilon}$, SCA Algorithm~\ref{alg_DCA} requires at most $\ceil{\bar{M}/\bar{\epsilon}}+1$ iterations. Empirically, however, our numerical experiments in Section~\ref{sec_Numerical_study} demonstrate that SCA Algorithm~\ref{alg_DCA} converges rapidly, often requiring far fewer iterations than this theoretical upper bound suggests.

\end{APPENDICES}

\end{document}